\documentclass[acmsmall, review=false]{acmart}

\makeatletter
\def\mdseries@tt{m}
\makeatother

\usepackage[utf8]{inputenc}
\usepackage[T1]{fontenc}
\usepackage{amsmath}
\usepackage{stmaryrd}
\usepackage{algorithm}
\usepackage[noend]{algpseudocode}
\usepackage{todonotes}

\usepackage[draft=true]{minted}

\usepackage{bm}

\usepackage{hyperref}
\title{Automatic Code Generation for High-Performance Discontinuous Galerkin Methods on Modern Architectures}
\author{Dominic Kempf}
\orcid{0000-0002-6140-2332}
\email{dominic.kempf@iwr.uni-heidelberg.de}
\author{René Heß}
\email{rene.hess@iwr.uni-heidelberg.de}
\author{Steffen Müthing}
\email{steffen.muething@iwr.uni-heidelberg.de}
\author{Peter Bastian}
\email{peter.bastian@iwr.uni-heidelberg.de}
\affiliation{
  \institution{Heidelberg University}
  \country{Germany}
}

\acmJournal{TOMS}

\newcommand\jump[1]{\llbracket #1 \rrbracket}
\newcommand\avg[1]{\{ #1 \}}
\newcommand\f[0]{\hspace{-0.1em}f\hspace{-0.15em}}
\renewcommand\vec[1]{\bm{#1}}

\begin{document}

\begin{abstract}
SIMD vectorization has lately become a key challenge in high-performance computing.
However, hand-written explicitly vectorized code often poses a threat to the software's sustainability.
In this publication we solve this sustainability and performance portability issue by
enriching the simulation framework dune-pdelab with a code generation
approach. The approach is based on the well-known domain-specific language UFL, but combines it
with loopy, a more powerful intermediate representation for the computational kernel.
Given this flexible tool, we present and implement a new class of vectorization strategies for
the assembly of Discontinuous Galerkin methods on hexahedral meshes exploiting the finite element's tensor product structure.
The optimal variant from this class is chosen by the code generator through an autotuning approach.
The implementation is done within the open source PDE software framework Dune and the discretization module dune-pdelab.
The strength of the proposed approach is illustrated with performance measurements for DG schemes for a scalar diffusion reaction equation
and the Stokes equation.
In our measurements, we utilize both the AVX2 and the AVX512 instruction set,
achieving 40\% to 60\% of the machine's theoretical peak performance for one matrix-free application of the operator.
\end{abstract}


\maketitle

\section{Introduction}

In the last decade, there have been two main approaches to finite element software frameworks in the academic community.
Library approaches like Dune \cite{bastian2008grid1} \cite{bastian2008grid2}, MFEM \cite{mfem} or deal.ii \cite{bangerth2007deal} use object-oriented C++ to provide users powerful toolboxes for numerical simulation codes.
The biggest advantage of this approach is its flexibility and extendability.
However, these frameworks have been struggling with the learning curve experience of new users, as well as the necessity for users to have quite some expertise in programming.
Code generation based frameworks such as FEniCS \cite{LoggMardalEtAl2012a} or firedrake \cite{rathgeber2015firedrake} pursued an alternate road focussing on ease-of-use:
The finite element problem is formulated in a domain specific language (DSL) and all user code is written through a high-level Python interface.
This provides an excellent framework for rapid prototyping, but often fails once the capabilities of the underlying C(++) framework need to be extended.
A related approach is implemented by projects like FreeFEM++, which embeds a DSL for PDEs into C++ using expression templates. \\

The necessity of expert level C++ programming skills has become even more apparent in the light of the latest hardware innovations.
Programming techniques that enable high performance computing on modern architectures include SIMD-aware programming and multi threading.
This work focusses on SIMD (Single Instruction Multiple Data) vectorization, which has become a crucial factor in achieving near-peak floating point performance on modern microarchitectures.
This is due to both the increase in bit width of SIMD registers and the ability to perform fused multiplication and addition in a single instruction.
As a consequence the maximum percentage of peak performance an optimal scalar code can achieve, has dropped significantly in the last decade:
While it could achieve 50\% of double precision peak performance on a 2008 Intel Nehalem processor (SSE4.2, no FMA), the number is as low as 6.125\% on the 2016 Intel Skylake processor family (AVX512, FMA). \\

The performance optimization literature describes two basic approaches to the issue (see e.g. \cite{franchetti2005simd}):
\begin{itemize}
 \item Automatic Vectorization through a vectorizing compiler. Developers write scalar code, whose data dependencies are analyzed by the compiler in order to identify vectorization opportunities.
 \item Explicit Vectorization, where developers directly write code for the target architecture.
\end{itemize}

While the former would be the most preferrable solution in terms of separation of concern between application developer and performance engineers, in many practical applications it is not feasible.
Many vectorization opportunities cannot be automatically found, as they would change the program semantics and only domain knowledge allows developers to include those into the optimization search space.
As an example we mention a memory layout change along the lines of Structure of Arrays vs. Array of Structures.
However, also explicit vectorization is suboptimal, as code is written directly for a given microarchitextures through low-level programming constructs such as compiler-specific intrinsics or even assembly routines.
The resulting code is usually hard to read, hard to maintain and hard to optimize further.
This poses a severe threat to the sustainability of the software.
We are targetting this issue by generating explicitly vectorized code from a hardware-independent problem formulation in a domain specific language. \\

Not all discretization methods for the solution of PDEs are equally well suited for high performance implementation on modern architectures.
In order to achieve a notable percentage of machine peak, the numerical algorithm needs to exhibit a favorable ratio of floating point operations per byte loaded from main memory.
We target high order Discontinuous Galerkin methods on hexahedral meshes for a variety of reasons beyond it being a powerful and flexible method for a large variety of PDE problems.
DG methods allow for globally blocked data structures. These data structures can be accessed directly from local computations, removing a costly data movement step from the algorithm.
With hexahedral meshes, the basis functions exhibit tensor product structure.
This can be exploited through sum factorization \cite{orszag1980spectral}, which not only greatly reduces the algorithmic complexity of the algorithm, but also features solely the fused multiplication and addition (FMA) operations needed to make full use of the floating point capabilities of modern CPUs.
The treatment of tensor product spaces in the domain specific language UFL has been described in \cite{mcrae2016tensorelements}.
Sum factorization has already been adopted into a variety of numerical codes e.g. \cite{muething2018sumfact} \cite{kronbichler2012genericinterface} \cite{schoeberl2014ngsolve} \cite{homoloya2017structure}.
Using sum facotrization and minimizing data movement, we can achieve a finite element assembly algorithm that is limited by the processors floating point capabilities instead of its memory bandwidth.
Replacing inherently memory-bound matrix-vector products with matrix free operator applications using the same sum factorization technique as in finite element assembly, we achieve compute boundedness even in the linear algebra part of our code.
In other words, our algorithms are designed to be capable of recomputing matrix entries faster than loading them from memory. \\

It is worth noting, that iterative solvers using this kind of matrix-free operator evaluation suffer from the additional challenge to implement preconditioners that do not hamper the algorithmic complexity of the overall algorithm.
Matrix free preconditioning techiques have been studied by various authors:
In \cite{mueller2018matrixfree}, Block-Jacobi and Block Gauss-Seidel preconditioners in a fully matrix-free and in a partially matrix-free setting are investigated.
\cite{pazner2018preconditioner} uses a Kronecker-product singular value decomposition approach to approximate the jacobian such that it can be evaluated matrix-free.
In \cite{diosady2017preconditioner}, the use of alternate-direction-implicit and fast diagonalization methods is advocated. \\

We are integrating the results of this publication into the Dune framework\cite{bastian2008grid1} \cite{bastian2008grid2}.
The Dune project emerged as a C++ successor to long standing C projects such as UG\cite{bastian1997ug}.
It features an large variety of grid managers available through a generic interface.
The discretization module dune-pdelab\cite{bastian2010pdelab} builds on top of the Dune framework and provides abstractions for finite element and finite volume discretizations.
Its strength lies with implicit methods and massive parallelism. \\

The contribution of this publication is twofold: A class of algorithms to achieve SIMD vectorization of the finite element assembly problem for DG methods with sum factorization is introduced.
These algorithms extend ideas presented in \cite{muething2018sumfact} to overcome their limitations with respect to the SIMD width and the applicable PDE problems.
Performance measurements for the Intel Haswell and Intel Skylake architectures show the algorithms cababilitity of achieving a notable percentage of peak performance.
The second contribution is the embedding of the algorithm into a code generation workflow to achieve performance portability across architectures and PDE problems.
This workflow is based on UFL, the state of the art domain specific language for the description of PDE models, but in contrast to the current workflow of other code generation based simulation packages uses a new intermediate representation.
This intermediate representation, loopy \cite{kloeckner2014loopy}, provides a powerful tool for optimizing the assembly loop nest and transforming otherwise hard-to-modify code properties like memory layouts.
The generated code is integrated into the simulation workflow of dune-pdelab through its CMake build system. \\

The structure of the paper is as follows:
The DG method and its implementation using sum factorization is summarized in section~\ref{sec:sumfact}.
Section~\ref{sec:vectorization} introduces new vectorization strategies extending on ideas from \cite{muething2018sumfact}.
These are implement through the code generation approach outlined in section~\ref{sec:codegen}.
The paper concludes with a performance study in section~\ref{sec:performance}.

\section{Sum Factorization for Discontinuous Galerkin Methods}
\label{sec:sumfact}

Section~\ref{sec:dg} will briefly introduce the notation and function spaces needed for Discontinuous Galerkin methods
with a special focus on tensor product structure of the finite elment. This will be done in a very general fashion, for practical examples
see section~\ref{sec:performance}. Section~\ref{sec:operatorappl} will then give a quick overview of how the tensor
product structure of the finite element can be exploited algorithmically through sum factorization.

\subsection{The Discontinuous Galerkin Method with Tensor Product Spaces}
\label{sec:dg}

Let $\{\mathcal{T}_h\}_{h>0}$ be a family of shape regular triangulations of the
domain $\Omega\subset \mathbb{R}^d$ consisting of closed elements $T$, each being
the image of a map $\mu_T : \hat T \to T$ with $\hat T$ being the reference
cube in $d$ dimensions. The map $\mu_T$ is differentiable, invertible and
its gradient is nonsingular on $\hat T$.
$F$ is an interior face if it is the intersection of two elements
$T^-(F), T^+(F)\in\mathcal{T}_h$ and $F$ has non-zero $(d-1)$-dimensional
measure.
All interior faces are collected in the set $\mathcal{F}_h^i$.
Likewise, $F$ is a boundary face if it is
the intersection of some $T^-(F)\in\mathcal{T}_h$ with $\partial\Omega$
and has non-zero $(d-1)$-dimensional measure.
All boundary faces make up the set $\mathcal{F}_h^{b}$
and we set $\mathcal{F}_h = \mathcal{F}_h^i \cup \mathcal{F}_h^{b}$.
The diameter of $F\in\mathcal{F}_h$ is $h_F$ and with each $F\in\mathcal{F}_h$
we associate a unit normal vector $\vec{n}_F$ oriented from $T^{-}(F)$ to $T^+(F)$
in case of interior faces
and coinciding with the unit outer normal of $\Omega$ in case of boundary faces.
Every face $F$ is the image of a map $\mu_F : \hat F \to F$ with $\hat F$ being the
reference element of the face. \\

Note that we restrict ourselves to cuboid meshes, where the reference element
is the tensor product cell $[0,1]^d$. We define the DG finite element space as the
tensor product of one dimensional polynomials on the reference interval, allowing
$v\in V_h^{\vec{k}}$ to be double-valued on $F\in\mathcal{F}_h^i$:
\begin{equation}
 V_h^{\vec{k}} = \left\{ v\in L^2(\Omega) : \forall\, T\in\mathcal{T}_h,
 v|_T = p \circ \mu_T^{-1} \text{ with } p\in\bigotimes_{i=0}^{d-1}\mathbb{P}^1_{k_i}([0,1])\right\}
\end{equation}
In this definition and throughout the rest of this work, the polynomial degree $\vec{k}=(k_0, k_1, \dots, k_{d-1})$
is allowed to be anisotropic, but we assume $\vec{k}$ to be constant throughout the triangulation.
This is not a general restriction of the presented methodology, we just have not tackled it yet
due to the involved programming efforts and the combinatorial increase in code size and generation time.
We define the number of local degrees of freedom per direction
as $n_i := \dim (\mathbb{P}^1_{k_i}([0,1])) = k_i + 1$ and the number of local degrees of freedom as
$N := \prod_{i=0}^{d-1}n_i$. For the presented work, the choice of the local basis of the space $\mathbb{P}^1_{k_i}([0,1])$
is not important.

We expect the discretized PDE problem to be given in residual form:
The solution $u_h\in V_h^{\vec{k}}$ is given as the solution of the discrete variational problem
\begin{equation}
 \label{equ:varprob}
 r(u_h, v_h) = 0 \ \ \forall v_h\in V_h^{\vec{k}}.
\end{equation}
We assume $r(u_h,v_h)$ to be expressable in the following form:
\begin{align}
\label{equ:splitteddg}
 r(u_h, v_h) & = \sum_{T\in\mathcal{T}_h}\int_{\hat{T}}\sum_{w\in\{v_h,\partial_0v_h,\dots ,\partial_{d-1}v_h\}} \tilde{r}_w^{volume}(\hat{u}_h, \mu_T, \hat{x})\ w\ d\hat{x}\\
 & + \sum_{F\in\mathcal{F}_h^i}\int_{\hat{F}}\sum_{\substack{w\in\{ v_h^+,\partial_0v_h^+,\dots ,\partial_{d-1}v_h^+,\\ v_h^-,\partial_0v_h^-,\dots ,\partial_{d-1}v_h^- \}}}\tilde{r}_w^{skeleton}(\hat{u}_h^+, \hat{u}_h^-, \mu_F, \hat{x})\ w\ d\hat{x} \\
 & + \sum_{F\in\mathcal{F}_h^b}\int_{\hat{F}}\sum_{w\in\{v_h,\partial_0v_h,\dots ,\partial_{d-1}v_h\}} \tilde{r}_w^{boundary}(\hat{u}_h, \mu_F, \hat{x})\ w\ d\hat{x}
\end{align}
This problem description introduces the notation that will be necessary in the following sections.
We do not expect the user to provide the input in this splitted form, but perform symbolic manipulation in our code generation toolchain to determine the functions $\tilde{r}$.
To solve problem~\eqref{equ:varprob} with Newton's method, we need not only be able to evaluate the algebraic residual $R_i(\vec{c}):=r(\sum_j c_j\phi_j ,\phi_i)$, but also its jacobian
$J_{ij}:=\frac{\partial R_i(\vec{c})}{\partial c_j}$. As mentioned, we would like to only implement the action $A_i(\vec{c}, \vec{z}) = (J\vec{z})_i$ of the jacobian on a vector $\vec{z}$.
Although we are omitting the exact formulae for brevity here, we mention that the action of the jacobian can be stated in a similar fashion as equation~\eqref{equ:splitteddg}.
For linear problems this will be the same as the residual $R_i(\vec{z})$ except for those terms not depending on $u_h$.
For nonlinear problems, the action of the jacobian will depend on both $\vec{c}$ and $\vec{z}$.

All integrals from equation~\eqref{equ:splitteddg} will be computed with appropriate quadrature rules.
We construct the quadrature rules on the reference cuboids by building the tensor product of 1D quadrature rules with maximum order.
The number of quadrature points per direction is defined as $m_i$ and is allowed to vary the same way as the polynomial degree: per direction, but not throughout the triangulation.
The total number of quadrature points on the reference cuboid is defined as $M:=\prod_{i=0}^{d-1}m_i$.

\subsection{Sum Factorized Residual Evaluation}
\label{sec:operatorappl}

In the following chapters we will focus on the evaluation of the residual $r(u_h,v_h)$ for the discussion of different vectorization strategies and performance results.
The described techniques are however fully applicable to matrix free jacobian applications and jacobian matrix assembly.
We will discuss the details of these generalizations at the appropriate places.
Note that, while we are formulating the algorithms here with a second order elliptic problem in mind, there is no hard assumptions on the residuum $\tilde{r}$ for the algorithms to be valid.
In particular, the residuum does not have to have tensor product structure. \\

The sum factorization technique is used in two places:
\begin{itemize}
 \item For the evaluation of finite element functions $\hat{u}_h$ on the reference element and its partial derivatives $\partial_i \hat{u}_h$ at all quadrature points.
 \item To multiply the evaluations of functions $\tilde{r}_w$ with $w$ (being the test function or its partial derivatives)
\end{itemize}

The necessary calculations are usually expressed in the literature using Kronecker products\cite{vanloan2000kronecker}.
We introduce the matrices $A^{(i)}\in\mathbb{R}^{m_i\times n_i}$, whose entries are the evaluations of the chosen local basis functions of $\mathbb{P}^1_{k_i}([0,1])$ at the given 1D quadrature points $\xi_0,\dots ,\xi_{m_i-1}$.
Given the coefficient vector $\vec{x}$ as a $d$-way tensor $X$, the evaluation of $\hat{u}_h$ at all quadrature points reads the following:
\begin{equation}
 \label{equ:withkronecker}
 \hat{U}=\left( A^{(0)}\otimes A^{(1)}\otimes\dots\otimes A^{(d-1)}\right) X
\end{equation}
For the evaluation of the partial derivative $\partial_j\hat{u}_h$, the basis evaluation matrix $A^{(j)}$ needs to be replaced with the matrix $D^{(j)}\in\mathbb{R}^{m_j\times n_j}$ containing the derivatives of the 1D basis functions at the 1D quadrature points.
We use the notation $\partial_j\hat{U}$ for the $d$-way tensor containing the $j$-th component of the gradient of $\hat{u}$ at all quadrature points.

The tensor $\hat{U}$ from equation~\eqref{equ:withkronecker} is evaluated in the following way, which is the fundamental idea of sum factorization:

\begin{align}
 \hat{U}_{i_0\dots i_{d-1}} &= \hat{u}_h(\bm{\xi}_{i_0\dots i_{d-1}}) \\
 &= \sum_{j_{d-1}=0}^{n_{d-1}-1}\cdots\sum_{j_1=0}^{n_1-1}\sum_{j_0=0}^{n_0-1}\prod_{k=0}^{d-1}A^{(k)}_{i_k,j_k}X_{j_0\dots j_{d-1}} \label{equ:highcomplex}\\
 &= \sum_{j_{d-1}=0}^{n_{d-1}-1}A^{(d-1)}_{i_{d-1},j_{d-1}}\cdots\sum_{j_1=0}^{n_1-1}A^{(1)}_{i_1,j_1}\sum_{j_0=0}^{n_0-1}A^{(0)}_{i_0,j_0}X_{j_0\dots j_{d-1}} \label{equ:lowcomplex}
\end{align}

Note how the complexity of the calculation reduces in equation~\eqref{equ:lowcomplex} in comparison to equation~\eqref{equ:highcomplex}.
Assuming $n_i=m_i=p$, the complexity of evaluating all elements of $\hat{U}$ decreases from $\mathcal{O}(p^{2d})$ to $\mathcal{O}(p^{d+1})$, illustrating well the desirability of the sum factorization approach.
On a linear algebra level, a sum factorization kernel boils down to a series of tensor contractions and tensor rotations.

The residual assembly process for one cell of a 3D problem is described in detail in algorithm~\ref{alg:assembly}.
Part 1 contains the aforementioned evaluation of finite element solution and its partial derivatives.
In part 2, the functions $\tilde{r}$ from equation~\eqref{equ:splitteddg} are evaluated at each quadrature point in order to set up an input tensor for part 3, which multiplies the test function in a sum factorized fashion.
To this end, transposed basis evaluation matrices $A^{(j), T}$ are necessary.

\begin{algorithm}
\begin{flushleft}
  \textbf{Input:} Tensor $X$ of local coefficients, $\mu_T$, Residual tensor $R$\\
  \textbf{Output:} Updated Residual tensor $R$
\end{flushleft}
 \begin{algorithmic}[1]
    \State $\hat{U}\gets\left( A^{(0)}\otimes A^{(1)}\otimes A^{(2)}\right) X$ \Comment{Part 1: Eval. of ansatz fct.}
    \State $\partial_0\hat{U}\gets\left( D^{(0)}\otimes A^{(1)}\otimes A^{(2)}\right) X$
    \State $\partial_1\hat{U}\gets\left( A^{(0)}\otimes D^{(1)}\otimes A^{(2)}\right) X$
    \State $\partial_2\hat{U}\gets\left( A^{(0)}\otimes A^{(1)}\otimes D^{(2)}\right) X$
    \For{$\hat{\xi}_{i_0i_1i_2} \in \text{quadrature points}$} \Comment{Part 2: Quadrature loop}
      \State $R^{v_h}_{i_0i_1i_2} \gets \tilde{r}_{v_h}^{volume}(\hat{U}_{i_0i_1i_2}, (\partial_0\hat{U})_{i_0i_1i_2}, (\partial_1\hat{U})_{i_0i_1i_2}, (\partial_2\hat{U})_{i_0i_1i_2}, \mu_T, \hat{\xi}_{i_0i_1i_2})$
      \State $(R^{\partial_0v_h})_{i_0i_1i_2} \gets \tilde{r}_{\partial_0v_h}^{volume}(\hat{U}_{i_0i_1i_2}, (\partial_0\hat{U})_{i_0i_1i_2}, (\partial_1\hat{U})_{i_0i_1i_2}, (\partial_2\hat{U})_{i_0i_1i_2}, \mu_T, \hat{\xi}_{i_0i_1i_2})$
      \State $(R^{\partial_1v_h})_{i_0i_1i_2} \gets \tilde{r}_{\partial_1v_h}^{volume}(\hat{U}_{i_0i_1i_2}, (\partial_0\hat{U})_{i_0i_1i_2}, (\partial_1\hat{U})_{i_0i_1i_2}, (\partial_2\hat{U})_{i_0i_1i_2}, \mu_T, \hat{\xi}_{i_0i_1i_2})$
      \State $(R^{\partial_2v_h})_{i_0i_1i_2} \gets \tilde{r}_{\partial_2v_h}^{volume}(\hat{U}_{i_0i_1i_2}, (\partial_0\hat{U})_{i_0i_1i_2}, (\partial_1\hat{U})_{i_0i_1i_2}, (\partial_2\hat{U})_{i_0i_1i_2}, \mu_T, \hat{\xi}_{i_0i_1i_2})$
    \EndFor
    \State $R\gets R+\left( A^{(0),T}\otimes A^{(1),T}\otimes A^{(2),T}\right)R^{v_h}$\Comment{Part 3: Mult. with test fct.}
    \State $R\gets R+\left( D^{(0),T}\otimes A^{(1),T}\otimes A^{(2),T}\right)R^{\partial_0v_h}$
    \State $R\gets R+\left( A^{(0),T}\otimes D^{(1),T}\otimes A^{(2),T}\right)R^{\partial_1v_h}$
    \State $R\gets R+\left( A^{(0),T}\otimes A^{(1),T}\otimes D^{(2),T}\right)R^{\partial_2v_h}$
 \end{algorithmic}

 \caption{Sum factorized algorithm for calculating the contributions of one 3D cell to the residual evaluation for a 2nd order elliptic problem.}
 \label{alg:assembly}
\end{algorithm}

Integrals over facets $F\in\mathcal{F}_h$ are treated similarly, only that $m_d:=1$ for $d$ being the normal direction on the reference cube.
It is worth to note, that in this case the total number of arithmetic operations of a sum factorization kernel can be greatly reduced by permuting the order of tensor contractions.
Such considerations of course make the implementation of a sum factorized algorithm dependent on the embeddings of $F$ into the reference elements of $T^-(F)$ and $T^+(F)$.
Code generation again helps greatly to write the implementations necessary for different face embeddings and combinations thereof.

\section{Explicit SIMD Vectorization}
\label{sec:vectorization}

In the following we will develop vectorization strategies for algorithm \ref{alg:assembly}.
These strategies follow the ideas previously described in \cite{muething2018sumfact}, but are extended to gain further flexibility w.r.t. new PDE models and meet the challenge of vectorizing for new microarchitectures.
We will classify our approaches into two categories: Loop-fusion based approaches and loop-splitting based ones.
Loop fusion based approaches typically require a drastical change in memory layout to work, where loop splitting based ones only work optimally if the mathematical problem leads to loop bounds satisfying suitable divisibility constraints.

\subsection{Loop Fusion Based Vectorization Strategies}
\label{sec:horizontal}

When identifying loops to fuse for a vectorization strategy in the finite element assembly problem, there are multiple choices for the level of granularity of the workload to be parallelized.
All these different approaches have advantages and disadvantages.
The most widely used approach (e.g. in \cite{kronbichler2012genericinterface}) is to parallelize over multiple cells in the grid.
The advantages of this strategy lie with the natural adoption to wider SIMD lanes on future architectures.
However, this also introduces additional costs:
The memory footprint of the integration kernel is increased by a factor of the size of SIMD lane width.
Additionally, the data for the degrees of freedom from multiple cells need to be interleaved in a preparation step.
This is an extension to a procedure that is commonly done in finite element assembly for continuous finite elements:
Degrees of freedom associated with a single element are gathered into a data structure for a local integration kernel to work on.
After the integration kernel has run, the data is scattered back into the global data structures.
However, Discontinuous Galerkin methods can also be implemented efficiently such that they operate on the blocks in the global data structure directly.
Doing so removes a large amount of costly, hard-to-hide memory operations from the algorithm and thus enables higher performance.
We have identified and described this issue in \cite{muething2018sumfact}.
Sticking to the idea of avoiding the setup of local data structures, we choose a different level of granularity to perform vectorization.
\\

Within an integration kernel on a cell or facet, usually several quantites need to be evaluated through a sum factorized algorithm.
Many of these sum factorization algorithms exhibit great structural similarities.
We explain the idea using the example of the residual evaluation algorithm~\ref{alg:assembly} for a second order elliptic PDE in 3D with a SIMD width of 256 bits.
We do this restriction to illustrate our core ideas and will then proceed how the approach generalizes to other models, architectures, space dimensions and jacobian assembly.
The core idea is to use the necessary four tensor quantities of step one in algorithm~\ref{alg:assembly} for vectorization:
The finite element solution $\hat{U}$ at all quadrature points $\xi_{i_0i_1i_2}$, as well as the three components of its gradient $\partial_k\hat{U}$.
We recall the main tensor product formulaes from section~\ref{sec:sumfact} for these quantities (in 3D):

\begin{align}
 \label{equ:unfusedtensor}
 \partial_0\hat{U} &= \left( D^{(0)}\otimes A^{(1)}\otimes A^{(2)}\right) X \\
 \partial_1\hat{U} &= \left( A^{(0)}\otimes D^{(1)}\otimes A^{(2)}\right) X \\
 \partial_2\hat{U} &= \left( A^{(0)}\otimes A^{(1)}\otimes D^{(2)}\right) X \\
           \hat{U} &= \left( A^{(0)}\otimes A^{(1)}\otimes A^{(2)}\right) X \label{equ:unfusedtensor_last}
\end{align}

As the tensor bounds $m_j$ and $n_j$ match in all of these computation, so do the loop bounds in the resulting sum factorization kernel implementation.
Therefore, the loops in the implementation can be fused to achieve an implementation suitable for SIMD vectorization.
In such implementation, each of the equations \eqref{equ:unfusedtensor} - \eqref{equ:unfusedtensor_last} would be carried out in one SIMD lane.
We will now introduce mathematical notation to reason about such a fused kernel as a tensor calculation.
To this end, we define two operations commonly used in tensor algebra \cite{vanloan2000kronecker}:
\begin{itemize}
 \item
The $vec$ operator maps a $d$-way tensor to a vector by flattening.
This operation imposes an order on the tensor axes.
This is completely analoguous to selecting strides in multi-dimensional array computation.
We will use this operator to refer to the representation of a tensor in memory.
 \item
Given $d$-way tensors $A_0,\dots ,A_n$ with identical bounds, we define $A_0|\dots |A_n$ as the $(d+1)$-way tensor constructed by stacking the tensors $A_i$ on top of each other.
We assume that the order of axes of the input tensors is preserved in the stacked tensor and that the new axis generated by stacking is the fastest varying.
Put in other words, the new axis has stride 1.
\end{itemize}

Using this notation, equations~\eqref{equ:unfusedtensor} to \eqref{equ:unfusedtensor_last} can be combined into one equation:
\begin{equation}
 \label{equ:fusedtensor}
 \partial_0\hat{U}|\partial_1\hat{U}|\partial_2\hat{U}|\hat{U} = \left( D^{(0)}|A^{(0)}|A^{(0)}|A^{(0)}\otimes A^{(1)}|D^{(1)}|A^{(1)}|A^{(1)}\otimes A^{(2)}|A^{(2)}|D^{(2)}|A^{(2)}\right)X|X|X|X
\end{equation}

\begin{figure}
 \centering
 \includegraphics[width=.6\textwidth]{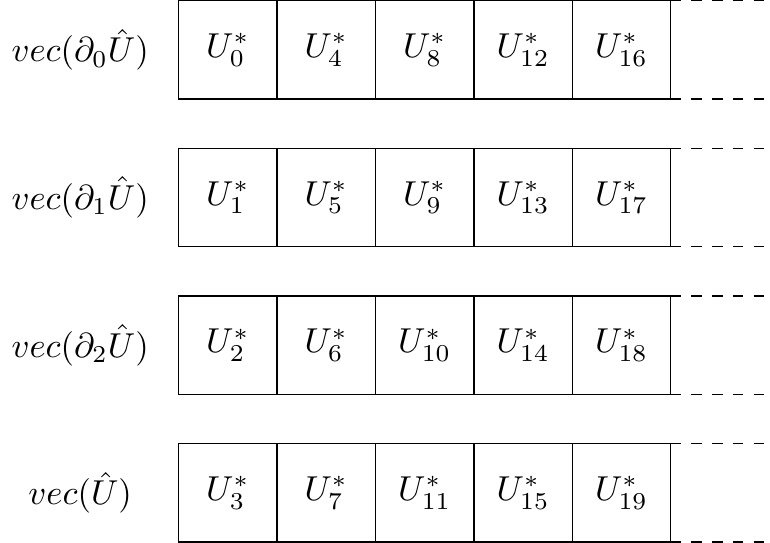}
 \caption{The memory layout of the tensor $U^*:=vec(\partial_0\hat{U}|\partial_1\hat{U}|\partial_2\hat{U}|\hat{U})$ is an array of structures obtained from the four original arrays $vec(\partial_0\hat{U})$, $vec(\partial_1\hat{U})$, $vec(\partial_2\hat{U})$ and $vec(\hat{U})$ by interleaving.}
 \label{fig:horizontal}
\end{figure}

We will now discuss the memory layout implications of implementing equation~\eqref{equ:fusedtensor}.
The memory layout of the input tensor $X$ is prescribed by the underlying discretization framework.
Accesses to the stacked tensor $X|X|X|X$ can be implemented easily by accessing an element of $X$ and broadcasting it into a SIMD register.
Our code generator will do this, whenever it finds a stacking of identical matrices.
The layout of the stacked basis evaluation matrices is given by interleaving the individual basis evaluation matrices, such that the stacked axis has stride $1$.
We do assemble these stacked basis evaluation matrices in memory.
This is a trade off decision between the increased memory traffic of loading redundant data and the necessity of instructions that manipulate single SIMD lanes.
These underlying matrices may be stored in column major or row major fashion, the code generation approach allows for flexibility in this regard.
Carrying out the sum factorization algorithm with these stacked tensors, the resulting stacked tensor $\partial_0\hat{U}|\partial_1\hat{U}|\partial_2\hat{U}|\hat{U}$ will be given as an interleaved tensor as well.
This can be seen as a an array of structures layout, where the inner structure is of fixed size $4$.
The layout is illustrated in figure~\ref{fig:horizontal}.
All data structures are aligned to the vector size to allow aligned loads into SIMD registers.
We will now discuss how this array of structures affects step 2 of algorithm~\ref{alg:assembly}, the quadrature loop.

\begin{figure}
 \centering
 \includegraphics[width=\textwidth]{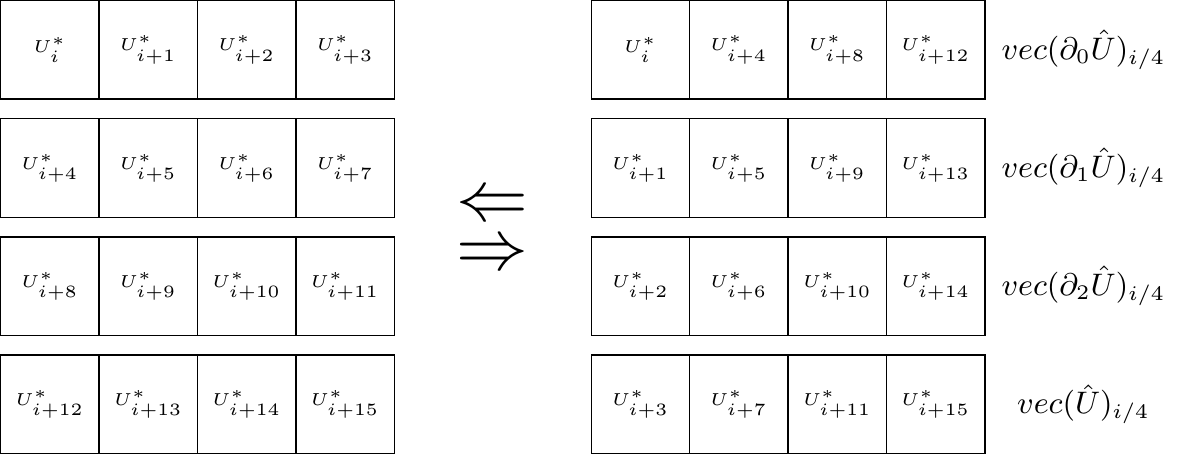}
 \caption{Register transposition needed in the quadrature loop: Four SIMD vectors of $U^*:=vec(\partial_0\hat{U}|\partial_1\hat{U}|\partial_2\hat{U}|\hat{U})$ are loaded and transposed in-place.
 The resulting four SIMD vectors have the layout needed for an efficient, straightforward vectorized implementation of the quadrature loop.
 The inverse operation is needed to get the correct layout for the input tensor for step 3 of algorithm~\ref{alg:assembly}.}
 \label{fig:transpose}
\end{figure}

Our idea of vectorizing the quadrature loop is based on the idea to treat four quadrature points at a time.
We have found it beneficial to neglect the tensor product structure of the quadrature loop here and instead use flat indexing.
In order to evaluate the functions $\tilde{r}$ from equation~\eqref{equ:splitteddg} with vector arguments, we need to undo
the array of structures layout. We do so in the quadrature loop by applying a transposition of four consecutive SIMD vectors
of $U^*:=vec(\partial_0\hat{U}|\partial_1\hat{U}|\partial_2\hat{U}|\hat{U})$.
The procedure is illustrated in figure~\ref{fig:transpose}.
The transposition code is implemented efficiently in C++ intrinsics. \\

Step three of algorithm~\ref{alg:assembly} again expects an array of structures type layout, that we do not get directly from evaluating the functions $\tilde{r}$ in a vectorized way.
We achieve this by applying the transposition algorithm from figure~\ref{fig:transpose} again.
The overall quadrature loop algorithm is summarized in algorithm~\ref{alg:quadloop}.
It is worth noting, that the $i$-loop does not need a tail loop although the total number of quadrature points is not necessarily divisible by four:
It is sufficient to overallocate the storage of $R$ to assure that the last loop iteration cannot write out of bounds.
Step three of algorithm~\ref{alg:assembly} is treated in the exact same way step one is, with the notable exception of the necessity to accumulate the results of four sum factorization kernel into the residual tensor.
With the chosen memory layout, this requires an intra-register reduction.
The implementation of this operation is subject to special care, as it benefits greatly from microarchitecture-dependent optimization.

\begin{algorithm}
 \caption{
Explicitly vectorized quadrature loop for a SIMD width of 4 and a total of $M$ quadrature points.
We use $U^*$ and $R^*$ as shortcuts for $vec(\partial_0\hat{U}|\partial_1\hat{U}|\partial_2\hat{U}|\hat{U})$ and
$vec(R^{\partial_0v}|R^{\partial_1v}|R^{\partial_2v}|R^v)$. The $j$ loop is implemented through SIMD
vectorization. This is formulated in the same frame as algorithm~\ref{alg:assembly}: Assembly of a 3D volume integral for a second order elliptic problem.
}
 \label{alg:quadloop}
 \begin{algorithmic}[1]
  \State $i\gets 0$
  \While{$i<M$}
    \State\Call{TransposeRegisters}{$U^*_{4i}$,\dots , $U^*_{4i+15}$}
    \For{$j\in\{ 0,1,2,3\}$}
      \State$R^*_{4i+j}\gets \tilde{r}^{volume}_{\partial_0v_h}(U^*_{4i+j}, U^*_{4(i+1)+j}, U^*_{4(i+2)+j}, U^*_{4(i+3)+j}, \mu_T, \hat{\xi})$
      \State$R^*_{4(i+1)+j}\gets \tilde{r}^{volume}_{\partial_1v_h}(U^*_{4i+j}, U^*_{4(i+1)+j}, U^*_{4(i+2)+j}, U^*_{4(i+3)+j}, \mu_T, \hat{\xi})$
      \State$R^*_{4(i+2)+j}\gets \tilde{r}^{volume}_{\partial_2v_h}(U^*_{4i+j}, U^*_{4(i+1)+j}, U^*_{4(i+2)+j}, U^*_{4(i+3)+j}, \mu_T, \hat{\xi})$
      \State$R^*_{4(i+3)+j}\gets \tilde{r}^{volume}_{v_h}(U^*_{4i+j}, U^*_{4(i+1)+j}, U^*_{4(i+2)+j}, U^*_{4(i+3)+j}, \mu_T, \hat{\xi})$
    \EndFor
    \State\Call{TransposeRegisters}{$R^*_{4i}$,\dots , $R^*_{4i+15}$}
    \State $i\gets i+4$
  \EndWhile
 \end{algorithmic}
\end{algorithm}

So far, we have studied the explicitly vectorized assembly algorithm under quite a number of simplifications and assumptions.
We will now discuss which of these assumptions were made for the sake of presentability in this article and which are actual limitations of the approach.
\begin{itemize}
 \item So far, we have only looked at volume integrals.
 Extension to boundary and interior facet integrals can be implemented straightforward by replacing the basis evaluation matrix for the normal direction of the facet
 with a special matrix consisisting of only one quadrature point. Note, that this quadrature point is either $0.0$ or $1.0$ depending on the facet being on the upper
 or lower boundary of the reference cube. Consequently, facet kernel implementations depend on the embedding of the facet into the reference cube.
 \item The approach was described for residual assembly only.
 Jacobian application for matrix-free methods works completely analoguous for linear problems.
 For nonlinear problems, not only the solution $\hat{u}_h$ needs to be evaluated in step 1 on of algorithm~\ref{alg:assembly}, but also the given linearization point.
 The evaluation of the linearization point uses the exact same algorithm and is vectorizable with the same techniques.
 Assembling a jacobian matrix is implemented in a very similar fashion:
 Step two and three of algorithm~\ref{alg:assembly} are executed once per basis function in the ansatz space.
 \item We have only looked at the weak formulation of an elliptic problem depending on both $\hat{u}_h$ and $\nabla \hat{u}_h$.
 For more complex problems, particularly systems of PDEs, the number of quantities that have to be evaluated via a sum factorization kernel may differ.
 This gives rise to additional opportunities to fuse sum factorization kernels.
 We offload the decision on how to vectorize to an algorithm implemented in our code generation workflow.
 The procedure is described in section~\ref{sec:sumfactcodegen}.
 \item We have only investigated the problem in three dimensional space.
 This is admittedly a best-case scenario, but we emphasize the relevance of such problems for real-world HPC applications.
 If the number of sum factorization kernels is less than what is needed for loop fusion (e.g. in 2D or because the evaluation $\hat{u}_h$ is not needed),
 it is still possible to perform above strategy with one SIMD lane left empty at a 25\% penalty to the floating point throughput.
 The following sections will introduce more clever strategies to cope with such situations though.
 \item We have limited ourself to a SIMD width of 256 bits.
 Instruction sets with smaller SIMD width like SSE2 can easily be targetted with the same techniques by combining four kernels into two vectorized kernels instead of one.
 Instruction sets with wider SIMD width like AVX512 pose a problem:
 Finding eight kernels in the problem that share the same loop bounds (or even the same input tensor) will rarely be possible.
 Developing vectorization techniques to overcome this limitation for AVX512 is the central goal of sections~\ref{sec:vertical} and \ref{sec:diagonal}.
 If these techniques are applied to another SIMD width $w$, a transposition algorithm for $w\times w$ matrices analoguous to figure~\ref{fig:transpose} has to be implemented.
 \item So far, we only fused sum factorization kernels, that share the same input tensor $X$, which resulted in values being broadcasted into SIMD registers.
 However, it is also possible to fuse kernels that have differing input tensors. This comes at the cost of an increased memory footprint of the fused kernel and more expensive
 load instructions (in the most general case, a gather instruction). In section~\ref{sec:sumfactcodegen} we will consider one special case that we consider a good trade-off between the increased costs and
 the desirable increase in vectorization opportunities: We fuse kernels with a total of two different input tensors.
 The necessary load instructions originate from the interoperability between different SIMD widths and fill the lower and upper half of a SIMD register.
 Use cases for having fusable kernels with exactly two input kernels occur quite naturally in DG discretizations:
 On a facet $F$, quantities need to be evaluated w.r.t. the cells $T^+(F)$ and $T^-(F)$ and when the action of the jacobian for a nonlinear problem is calculated, both the finite element solution and the linearization point need to be evaluated.

\end{itemize}

\subsection{Loop Splitting Based Vectorization Strategies}
\label{sec:vertical}

While the loop fusion vectorization from section~\ref{sec:horizontal} tries to fuse multiple sum factorization kernel, the idea of this section is to split the workload of one sum factorization kernel, such that execution can make use of SIMD parallelism.
This does not suffer the disadvantages seen above: Increased memory footprint of the kernel and the necessity of memory layout adjustments.
On the other hand, these splitting based vectorization techniques come with the disadvantage, that maximum efficiency can only be reached if the kernel structure exhibits loop bounds with suitable divisibility constraints.
We will now explore these strategies.
For the sake of readability, we again limit ourselves to a SIMD width of four lanes and to the evaluation treatment of the evaluation of $\hat{U}$ via
\begin{equation}
 \label{equ:evalprob}
 \hat{U}=\left( A^{(0)}\otimes A^{(1)}\otimes A^{(2)}\right) X
\end{equation}
and discuss possible generalizations later. \\

\begin{figure}
  \centering
  \includegraphics[width=.7\textwidth]{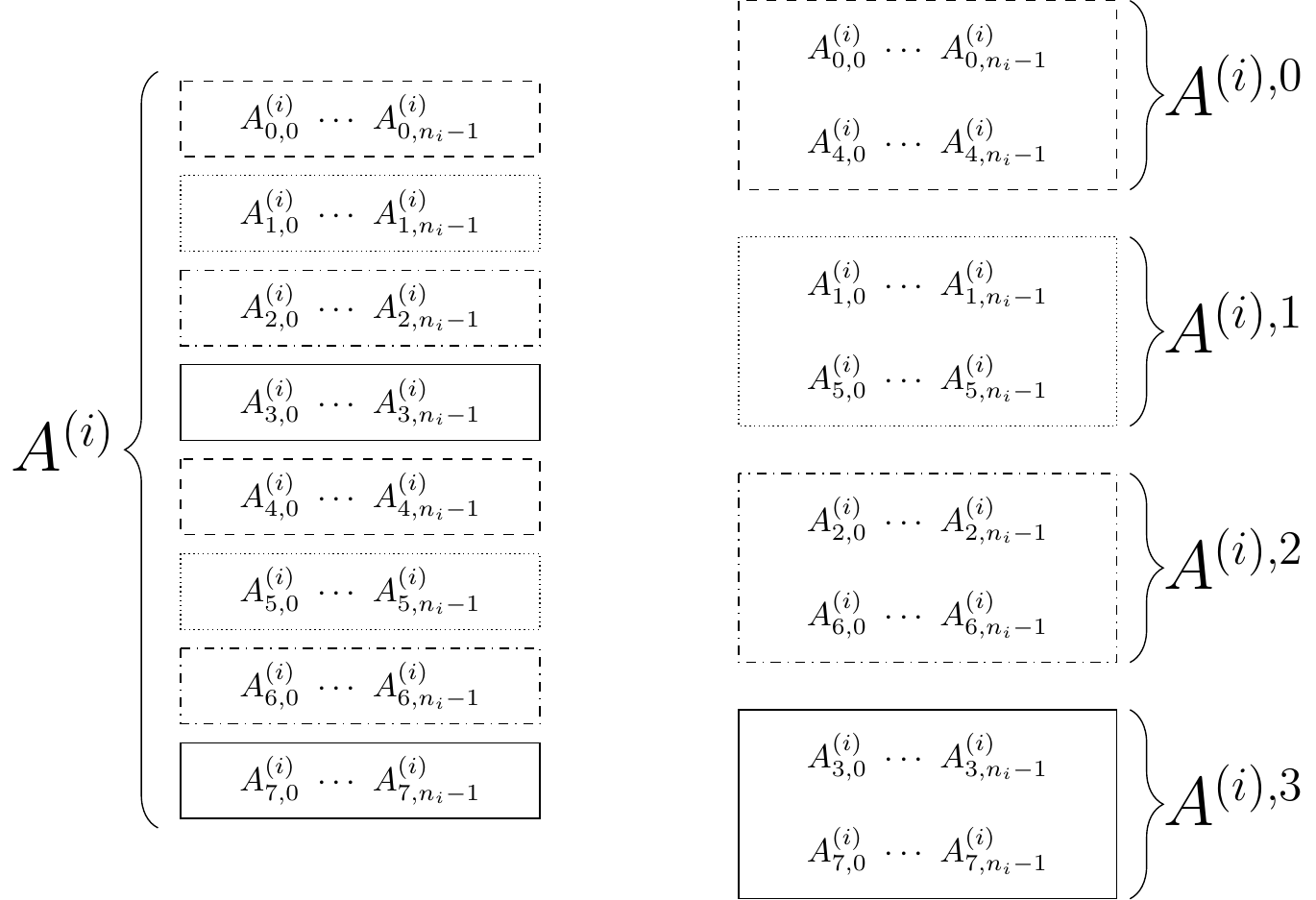}
  \caption{Slicing of the basis evaluation matrix: $A^{i}$ is split into four matrices $A^{(i), s}$ with $s\in\{ 0,\dots ,w-1\}$ in a circular fashion.}
  \label{fig:slicing}
\end{figure}

We base our strategy on the idea to split the set of quadrature points into a number of subsets equal to the SIMD width.
We do so by choosing one direction $i$ and splitting the basis evaluation matrix $A^{(i)}$ into $w$ matrices, where $w$ is the SIMD width.
For now, we assume the number of quadrature points to be divisible by $w$ and discuss other cases later on.
The index $s\in\{ 0,\dots ,w-1\}$ denotes the index of the slice $A^{(i), s}$ of the basis evaluation matrix.
Note, that we did not split $A^{(i)}$ in a blocked fashion, but in a circular one, like it is illustrated in figure~\ref{fig:slicing}.
Carrying out the sum factorized computation from equation~\eqref{equ:evalprob} with a slice $A^{(i), s}$ instead of $A^{(i)}$ leads to an output tensor $\hat{U}^s$ only containing evaluations at $1/w$ of the quadrature points.
Problem~\eqref{equ:evalprob} can then be reformulated into the following equivalent problem using the notation from section~\ref{sec:horizontal}:
\begin{equation}
 \label{equ:evalvec}
 \hat{U}^0|\hat{U}^1|\hat{U}^2|\hat{U}^3 = \left( A^{(0),0}|A^{(0),1}|A^{(0),2}|A^{(0),3}\otimes A^{(1)}|A^{(1)}|A^{(1)}|A^{(1)}\otimes A^{(2)}|A^{(2)}|A^{(2)}|A^{(2)}\right) X|X|X|X
\end{equation}
The fact that we have sliced $A^{(0)}$ in a circular fashion leads to the following, desirable property that allows us to load data of the resulting tensor $\hat{U}^0|\hat{U}^1|\hat{U}^2|\hat{U}^3$ without further manipulation of memory layout:
\begin{equation}
 vec(\hat{U}^0|\hat{U}^1|\hat{U}^2|\hat{U}^3)=vec(\hat{U})
\end{equation}

We observe that in contrast to section~\ref{sec:horizontal}, the combined basis evaluation matrices do not have to be explicitly set up beforehand, as $vec(A^{(0),0}|A^{(0),1}|A^{(0),2}|A^{(0),3})=vec(A^{(0)})$ and $A^{(i)}|A^{(i)}|A^{(i)}|A^{(i)}$ can be implemented as a broadcast of elements of $A^{(i)}$.
The loop splitting based approach described in this section does not depend on the problem structure the same way as the loop fusion based one from section~\ref{sec:horizontal} does.
As there is no need to group multiple sum factorization kernels, the approach vectorizes equally well in two and three dimensional space, as well as with arbitrary combination of terms present in the problem formulation.
However, applicability of the approach depends on the divisibility of the number of 1D quadrature points.
Having this constraint on the number of quadrature points is not as bad as having it on the number of basis functions:
Artificially increasing the number of quadrature points is equivalent to overintegration, which even yields additional accuracy for problems that are not exactly integrated.
However, one has to be cautious as the increase in floating point operations affects the whole algorithm and not only the sum factorization kernel to be vectorized.
We now study the cost increase of such a procedure and refer to section~\ref{sec:sumfactcodegen} for discussion of the necessary trade off decisions. \\

We recall that the number of quadrature points per direction is given as a tuple $\mathbf{m}=(m_0,\dots ,m_{d-1})$ and the number of basis functions per direction as a tuple $\mathbf{n}=(n_0,\dots ,n_{d-1})$.
The floating point cost $\mathcal{C}^{SF}(\mathbf{m}, \mathbf{n})$ of a single sum factorization kernel reads

\begin{equation}
 \label{equ:sumfactcost}
 \mathcal{C}^{SF}(\mathbf{m}, \mathbf{n}) = 2\sum_{k=0}^{d-1}\prod_{i=0}^km_i\prod_{j=k}^{d-1}n_j.
\end{equation}

We observe that $\mathcal{C}^{SF}(\mathbf{m}, \mathbf{n})$ is linear in $m_0$.
This also holds for other relevant parts of the algorithm such as the quadrature loop.
Consequently, the total cost of the algorithm will be increased by a factor of $\frac{m^*}{m_0}$, if the number of quadrature points in the first direction is increased to $m^*$.
Setting $m^*$ to the next multiple of the SIMD width $w$, we get a cost increase of $\lceil \frac{m_0}{w} \rceil \cdot \frac{w}{m_0}$.
For sufficiently high numbers of quadrature points, this increase becomes negligible.
In the worst case scenario of a $\mathbb{Q}_1$ discretization with minimal quadrature order, it can be as high as $\frac{w}{2}$ though and careful consideration is necessary.
See section~\ref{sec:sumfactcodegen} for more details about that. \\

In this section, we formulated the implementation idea in terms of the fusion based vectorization described in section~\ref{sec:horizontal}.
The same ideas could have been developed from a different perspective, but having it formulated using the same notation will be of great benefit for developing hybrid strategies in section~\ref{sec:diagonal} and also for the vectorization heuristics in the code generator described in section~\ref{sec:sumfactcodegen}. \\

\subsection{Hybrid Strategies}
\label{sec:diagonal}

Neither the strategy described in section~\ref{sec:horizontal} nor the strategy from section~\ref{sec:vertical} are in general a perfect fit for vectorization with wider SIMD width.
For the loop fusion strategy, the problem description will usually not exhibit enough quantities that can be computed in parallel.
The loop splitting strategy on the other hand leads to prohibitively severe constraints on the number of quadrature points with increasing SIMD width.
We now seek to combine the two strategies into a hybrid vectorization strategy mitigating these disadvantages. \\

We have formulated the loop fusion approach from section~\ref{sec:horizontal} and the loop splitting one from section~\ref{sec:vertical} using a common framework:
A set of sum factorization kernels is collected into a larger tensor, which is suitable for vectorization, where these kernels are potentially gained by first splitting the given sum factorization kernels.
We will now generalize this for arbitrary SIMD width and combinations of these techniques.
We define $f$ and $s$, such that $f$ denotes the number of sum factorization kernels to be combined through loop fusion and $s$ denotes the number of slices these are split into.
We only treat those cases, where $f\cdot s=w$ with $w$ the number of SIMD lanes.
For $f=4$ and $s=2$, this allows a natural extension of section~\ref{sec:horizontal} for a SIMD width $w=8$ (AVX512), which calculates $\hat{u}_h$ and $\nabla \hat{u}_h$ in parallel and introduces a divisibility constraint of $2$ on $m_0$:

\begin{align}
 \partial_0\hat{U}^0|\partial_0\hat{U}^1|\partial_1\hat{U}^0|\partial_1\hat{U}^1| & \partial_2\hat{U}^0|\partial_2\hat{U}^1|\hat{U}^0|\hat{U}^1 \nonumber\\
   = ( & D^{(0), 0}|D^{(0), 1}|A^{(0), 0}|A^{(0), 1}|A^{(0), 0}|A^{(0), 1}|A^{(0), 0}|A^{(0), 1} \nonumber\\
 & \otimes A^{(1)}|A^{(1)}|D^{(1)}|D^{(1)}|A^{(1)}|A^{(1)}|A^{(1)}|A^{(1)} \nonumber\\
 & \otimes A^{(2)}|A^{(2)}|A^{(2)}|A^{(2)}|D^{(2)}|D^{(2)}|A^{(2)}|A^{(2)}) \nonumber\\
 &X|X|X|X|X|X|X|X \label{equ:hybrid}
\end{align}

Again, we will study the memory layout implications of implementing equation~\eqref{equ:hybrid}.
The stacked basis evaluation matrices are preevaluated and loaded from memory, just like in section~\ref{sec:horizontal}.
The input tensor $X|X|X|X|X|X|X|X$ is implemented by broadcasting the values of $X$.
The only remaining question is how the memory layout of the output tensor $\partial_0\hat{U}^0|\partial_0\hat{U}^1|\partial_1\hat{U}^0|\partial_1\hat{U}^1|\partial_2\hat{U}^0|\partial_2\hat{U}^1|\hat{U}^0|\hat{U}^1$ affects the quadrature loop implementation.
Independently of the choice of $f$ and $s$, the quadrature loop treats $w$ quadrature points at a time.
However, to get $w$ values of the $f$ quantities present in the data, we need to shuffle $f$ consecutive vectors.
This results in the need for non-square matrix shufflings to get the correct data layout.
In order to reduce the amount of necessary such transposition implementations, we fix the intra-register layout to be such, that kernels resulting from splitting need to be adjacent to each other.
In other words, we disallow tensors like $\partial_0\hat{U}^0|\partial_1\hat{U}^0|\partial_2\hat{U}^0|\hat{U}^0|\partial_0\hat{U}^1|\partial_1\hat{U}^1|\partial_2\hat{U}^1|\hat{U}^1$.
We implement all necessary shuffling operations in C++ intrinsics.
The quadrature loop algorithm~\ref{alg:hybridquadloop} is further complicated by the fact, that an integration kernel might consist of more than one vectorized sum factorization kernel and that the choice of $f$ and $s$ can differ for each of these.

\begin{algorithm}
  \caption{General quadrature loop with a hybrid vectorization strategy for SIMD width of $w$:
    The input data is given as set of flat tensors $\{U^*_j\}$. Similarly, the output will be written as a set of flat tensors $\{R^*_k\}$.
    Each of of these flat tensors results from a fusion of a set of $w$ sum factorization kernels.
    The function $f$ maps those tensors to the number of fused quantities as described in section~\ref{sec:diagonal}.}
 \label{alg:hybridquadloop}

\begin{algorithmic}[1]
  \State $i\gets 0$
  \While{$i < M$}
    \For{$j\in\{ 0,\dots\}$}
      \State\Call{TransposeRegisters}{$(U^*_j)_{f(U^*_j)i}, \dots, (U^*_j)_{f(U^*_j)i + f(U^*_j)w -1}$}
    \EndFor
    \For{$k\in\{ 0,\dots\}$}
      \State$R^*_k\gets\text{Quadrature computation}$
      \State\Call{TransposeRegisters}{$(R^*_k)_{f(R^*_k)i},\dots , (R^*_k)_{f(R^*_k)i + f(R^*_k)w - 1}$}
    \EndFor
    \State $i\gets i+w$
  \EndWhile
 \end{algorithmic}
\end{algorithm}

We have seen, that the techniques of sections~\ref{sec:horizontal} and \ref{sec:vertical} can be combined to mitigate their disadvantages and target wider SIMD widths.
However, for a given problem, the number of possible vectorization strategies is of exponential nature and it is not a priori known, which one is best.
This issue is targetted in section~\ref{sec:sumfactcodegen} with a cost model approach.

\section{Performance Portability Through Code Generation}
\label{sec:codegen}

After studying vectorization techniques for variable SIMD widths in section~\ref{sec:vectorization},
we will now look into how the sustainability issue arising from the inherent hardware dependency
can be solved. To that end, we will first establish a code generation workflow in section~\ref{sec:toolchain}.
Then, we will see how optimal vectorization can be found in this workflow in section~\ref{sec:sumfactcodegen}.

\subsection{Code Generation Tool Chain}
\label{sec:toolchain}

\begin{figure}
 \centering
 \includegraphics[width=.7\textwidth]{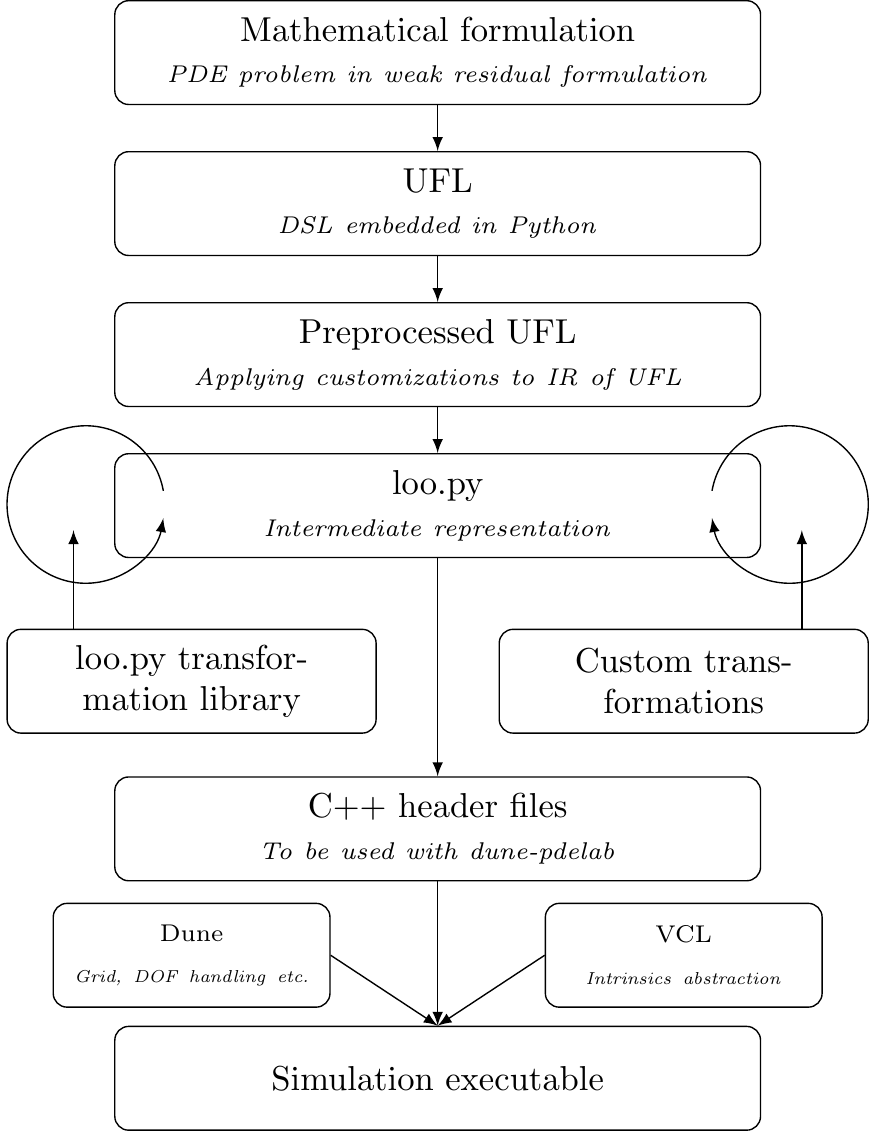}
 \caption{The employed code generation tool chain at a glance}
 \label{fig:toolchain}
\end{figure}

We are targetting dune-pdelab with our code generation workflow.
It builds upon the Dune core modules, which provide Dune's generic grid interface\cite{bastian2008grid2},
linear algebra, local basis functions, quadrature rules and geometry mappings.
The discretization framework dune-pdelab conceptionally extends this with key ingredients of finite element programs:
\begin{itemize}
 \item A \textit{FiniteElement} describes a local finite element on a given reference element.
 This comprises the local basis and how the local coefficients are attached to subentities of the reference element.
 \item A \textit{FiniteElementMap} maps grid cells to their associated local finite elements.
 One of PDELab's strengths is to exploit this knowledge at compile time and only pay runtime penalties if necessary (e.g. in $hp$-adaptivity).
 \item A \textit{GridFunctionSpace} implements a discrete ansatz space by combining a grid with a finite element map.
 PDELab allows to construct arbitrary trees of grid function spaces to represent systems of PDEs.
 Again, PDELab draws its strength from compile time reasoning about these tree structures allowing
 a lot of different blocking and ordering techniques of the global data structures.
 \item A \textit{ConstraintsAssembler} implements how degrees of freedom need to be constrained depending on boundary conditions and parallelism concepts.
 \item The finite element assembly algorithm is driven by a \textit{GridOperator} that iterates over the grid elements and facets.
 Implementation of assembly integrals is delegated to a \textit{LocalOperator} that only operates on single cells/facets.
 The grid operator also makes sure to apply constraints as needed.
\end{itemize}
With the intended code generation workflow, we would like to leverage as many features of the underlying framework as possible.
We therefore concentrate on generating code for the local operator, which is the performance critical component for finite element assembly.
The generation and compilation process is controlled by an extension to Dune's CMake build system.
We provide rudimentary automation of the full simulation workflow only for the very important topic of automated testing \cite{kempf2015testtools}, but advise users to still write their own Dune applications.
Otherwise, we'd have to make sure to deliver Python bindings to all of the frameworks features or we defeat the purpose of a highly modular C++ framework.

In the following we describe the employed tool chain in detail, which is also illustrated by figure~\ref{fig:toolchain}.
We aim to reuse existing code generation projects wherever feasible.
The FEniCS project has developed UFL \cite{AlnaesEtAl2012}, a domain specific language for the description of the weak formulations arising in finite element discretizations.
It is embedded in Python and also used outside of FEniCS, e.g. in the firedrake project.
Reusing UFL as the input language to our code generation toolchain will not only save us work, but also contribute to a standardization of open source finite element packages.
UFL allows to express multilinear forms in a cell- and facet-local fashion.
The Discontinuous Galerkin methods we are targetting in this work fit nicely into this framework. \\

UFL delivers an intermediate representation (IR) of the PDE problem in form of an abstract syntax tree (AST).
We apply some customizations and algorithms to the IR of UFL as a preprocessing:
\begin{itemize}
  \item
We enforce that users write the weak formulation in residual form (which is always a rank one form), as dune-pdelab's abstractions are formulated in this way.
This is quite similar to the way that users implement non-linear problems in FEniCS, though we use the same workflow for linear problems as well
\footnote{For readers that are more involved into UFL: We enforce this by overriding the \mintinline{python}{TrialFunction} class from UFL to be a \mintinline{python}{Coefficient} of reserved index, instead of an \mintinline{python}{Argument}.}.
From this residual form, we derive the bilinear form in the code generation process by taking the Gâteaux derivative with respect to the trial function with the automatic differentiation code provided by UFL.
For matrix free calculations, we apply another symbolic manipulation algorithm from UFL to get a rank one form implementing the action of the bilinear form on a given vector.
These symbolic manipulations do add a lot of value to our user experience, as the problem needs to be implemented exactly once and the jacobians can be derived from that.
In handwritten dune-pdelab codes, one can either pay a performance penalty and rely on automatic differentiation of residuals at runtime or one needs to handcode jacobians and their action as well.
Especially for nonlinear problems, the latter can be quite tedious.
  \item
We manipulate the given form such that we get symbolic representations of the functions $\tilde{r}^{volume}_w$ from equation~\eqref{equ:splitteddg}, where $w\in W:=\{ v_h,\partial_0v_h,\dots ,\partial_{d-1}v_h\}$.
A function $\tilde{r}^{volume}_w$ is extracted by traversing the form and replacing $w$ with $1$ and all other elements of $W$ with $0$.
Skeleton and boundary terms work the same with $W$ being defined as in equation~\eqref{equ:splitteddg}.
\item
We apply knowledge about geometry transformations for axiparallel grids.
These are not built into UFL, as the underlying discretization framework of FEniCS did not support hexahedral meshes for the longest time.
\end{itemize}

The IR of UFL is not a good choice for employing transformation-based optimization.
This mainly stems from the fact, that the only notion of loops in the UFL AST are sum reductions.
However, hardware-dependent transformations need much more insight into the loop structure of the assembly kernel.
We therefore transform the UFL IR into loopy \cite{kloeckner2014loopy} adding the necessary loop bound information for quadrature rules, ansatz and test functions.
Loopy is a Python package which provides an IR for computational kernels and a transformation library to operate on that IR.
The IR comprises a polyhedral description of the loop domain and a symbolic representation of the calculation to be carried out in the form of statements.
The package has already been proven to be capable of handling examples as complicated as full PDE applications \cite{kloeckner2016femexample}. \\

Loopy provides several code generation backends called targets, such as plain C, Cuda, OpenCL, ISPC etc.
The object-oriented nature of the code generation target classes allows to customize a target for C/C++ to instead produce code to be used with the discretization framework dune-pdelab without modifying loopy itself.
We use this to produce C code for our finite element assembly kernels that features enough C++ to call framework functions.
The loopy target for plain C code does not support explicit vectorization, as SIMD vectors are not part of the C language.
We add this functionality by generating code, that uses a C++ vectorization library.
We chose the vector class library \cite{fogvcl}, although there are several other viable options (e.g. VC \cite{kretz2012vc}).
The library provides C++ classes for a given precision and vector width (like \mintinline{c++}{Vec4d} for 256-bit double precision) representing a vector register having suitable operator and function overloads for all the basic tasks.
The implementation uses C++ intrinsics and hides most of the hardware-specific assembly-level considerations from the user and our code generator.
We consider such a library a key ingredient to a sustainable, performance portable code generation workflow. \\

Having settled on using loopy as the working horse of our code generation approach, we also gain additional value through loopys transformation library.
Common performance optimization techniques such as loop tiling, loop fusion or software prefetching are readily available.
We stress that in order to be capable of implementing such transformations, an IR needs to operate on a very specific level of abstraction:
On the one hand it needs to be reasonably general in the sense that it has a full symbolic representation of the computations to be carried out.
On the other hand it needs to have a notion of loop domains without fixing the loop structure too early in the tool chain (like an AST for the C language would).
We believe that loopy's level of abstraction hits this sweet spot needed for a high performance code generation workflow. \\

Our approach differs from other approaches to code generation in UFL-based projects:
The FEniCS project generates C code through the form compiler FFC\cite{kirby2006ffc} and any performance considerations are left to the underlying C++ framework dolfin and the C compiler.
Also, FEniCS - successfully - emphasizes usability and productivity over performance, which does not match with our intent of using code generation.
The firedrake project uses the form compiler TSFC \cite{homoloya2017tsfc}:
In a first step, UFL input is transformed into the tensor algebra IR GEM, which does not contain any domain-specific finite element information.
Such a step also happens in our transition from UFL to loopy, as loopy represents tensor algebra through the package pymbolic.
In a second step, TSFC generates C code from the given tensor algebra expression, much like the scheduler of loopy does.
However, this step does not yet take into account any hardware-specific considerations and performance optimization is left to COFFEE\cite{luporini2015coffee}, which operates on a C-style IR.
This toolchain is strictly following the idea of separation of concern, which we believe to not be well-applicable for performance optimization of finite element applications.
The firedrake project has recently also invested much work into generating code for sum factorized finite element assembly \cite{homoloya2017structure}. \\

Our code is published as a Dune module under a BSD license at
\begin{center}
 \url{https://gitlab.dune-project.org/extensions/dune-codegen}.
\end{center}

\subsection{Generating Explicitly Vectorized Sum Factorization Kernels}
\label{sec:sumfactcodegen}

We will now take a closer look at how the vectorization strategies described in section~\ref{sec:vectorization} can be exploited from the code generation tool chain.
Foremost, this is about how to choose the best vectorization strategy from a large set of possibilities.
All of the AST transformations (UFL to loopy, loopy to C) are implemented via recursive tree traversals with type-based function dispatch.
These algorithms work best if the tree transformation is fully local, meaning that the visitor object is completely stateless.
As our vectorization strategies - especially those based on loop fusion - are inherently non-local, we instead apply a two-step procedure:
\begin{itemize}
 \item During tree traversal, any quantity that is calculated through a sum factorization kernel is represented by a dedicated AST node \mintinline{python}{SumfactKernel}, that stores all the relevant information.
 \item After tree traversal, an algorithm decides on vectorization by providing a mapping of all the \mintinline{python}{SumfactKernel} in the AST to \mintinline{python}{VectorizedSumfactKernel} nodes with vectorization information attached.
 \item A second tree traversal is done, in which these modified \mintinline{python}{VectorizedSumfactKernel} nodes are realized by loopy statements.
\end{itemize}

We will now give details about the algorithm used to decide which vectorization strategy should be in use.
The number of possibilities to combine the vectorization strategies described in section~\ref{sec:vectorization} for a given set of sum factorization kernels is vast and it is not a priori known which strategy delivers optimal performance.
We mention two scenarios arising in the examples in section~\ref{sec:performance} to illustrate that trade off decisions need to be made:
\begin{itemize}
 \item
For a simple Poisson problem in 3D, $\partial_iu$ are needed, but not the evaluation of $u$ itself.
Given a SIMD width of 256 bits, is it better to fuse three kernels and ignore the forth SIMD lane or to apply a splitting based vectorization strategy partially or fully?
How does the quadrature order affect this?
 \item
For the implementation of the $\mathbb{Q}_k/\mathbb{Q}_{k-1}$ DG scheme for the Stokes equation from section~\ref{sec:stokes}, the evaluation of pressure cannot be parallelized with any other necessary evaluations.
Vectorizing pressure evaluation by splitting may come at the cost of increasing the number of quadrature points for the whole algorithm though.
When is it better to not vectorize pressure evaluation?
\end{itemize}
A costmodel based approach is required in order to make optimal vectorization decisions.
Such an approach consists of two core components: A function that systematically traverses all vectorization opportunities to find a minimum and an actual cost function.
In order to handle the exponential complexity of the traversal of vectorization opportunities, we employ a divide and conquer strategy splitting the optimization problem into several subproblems:
\begin{itemize}
 \item
Starting from the quadrature point tuple $(m_0,\dots ,m_{d-1})$ that was specified by the user or deduced from the problem formulation, we list all possible tuples with increased number of quadrature points, that enable other vectorization strategies.
For each of those we find an optimal vectorization strategy and find the minimum among these:
\begin{equation}
\label{equ:optimize_quad}
\begin{aligned}
& \underset{q}{\text{minimize}}
& & \text{cost}(\text{FixedQPMinimalStrategy}(\text{sumfacts}, \text{width}, q)) \\
&
& & q\in\left\{\left(\left\lceil\frac{m_0}{i}\right\rceil \cdot i, \dots, m_{d-1}\right)\ \middle|\ i=1,2,4,\dots, w\right\}.
\end{aligned}
\end{equation}
Here, $w$ again denotes the SIMD width.
\item
When finding an optimal strategy for a given fixed quadrature point tuple, we first divide the given set of sum factorization kernels into smaller subsets, which may potentially be subject to a loop fusion based vectorization approach i.e. they share the same loop bounds.
Minimal solutions w.r.t. the defined cost function for these subsets are then combined into a full vectorization strategy.
\end{itemize}

We define a function $parallelizable(s\f{})$ such that two sum factorization kernels that yield the same value are potentially vectorizable via loop fusion.
Likewise, we define a function $input(s\f{})$ such that two kernels yielding the same result operate on the same input tensor.
To wrap up the divide and conquer approach, a function $combine$ that merges the minimal solutions on subsets is used.
Algorithm~\ref{alg:vecopt} shows the overall optimization algorithm, where algorithm~\ref{alg:vecgen} shows the algorithm within one of the subsets.

\begin{algorithm}
 \caption{
Finding an optimal vectorization strategy: Given a set of sum factorization kernels, the SIMD width and a fixed quadrature point tuple, a vectorization strategy is returned as a mapping of sum factorization kernels to vectorized kernels.}
 \label{alg:vecopt}

 \begin{algorithmic}[1]
  \Function{FixedQPMinimalStrategy}{$sum\f{}acts$, $width$, $q$}
    \State $groups\gets \emptyset$
    \ForAll{$s\f{} \in sum\f{}acts$}
      \State insert $s\f{}$ in $groups[parallelizable(s\f{})]$)]
    \EndFor
    \State $results\gets\emptyset$
    \ForAll{$groupsum\f{}acts \in \text{values}(groups)$}
      \State $vecs\f{}\gets$ \Call{VectorizationStrategies}{$groupsum\f{}acts$, $width$, $q$}
      \State insert $\text{argmin}_{v\in vecs\f{}}\text{cost}(v)$ in $results$
    \EndFor
    \State\Return combine($results$)

  \EndFunction
 \end{algorithmic}
\end{algorithm}

\begin{algorithm}
 \caption{For a given set of sum factorization kernels, that are pairwise implementable in parallel, return all vectorization strategies from the pool of implemented methods in section~\ref{sec:vectorization}.
 The combine function merges the given mappings into one large mapping.}
 \label{alg:vecgen}

 \begin{algorithmic}[1]

 \Function{VectorizationStrategies}{$sum\f{}acts$, $width$, $q$}
   \State $strategies\gets\emptyset$

   \State $input\_groups\gets\emptyset$
   \ForAll{$s\f{}\in sum\f{}acts$}
     \State insert $s\f{}$ in $input\_groups[input(s\f{})]$
   \EndFor

   \ForAll{$num\_inputs\in \{1,2\}$}
     \If{$num\_inputs > len(input\_groups)$}
       \State break
     \EndIf

     \ForAll{$\f{}\in\{ 1,2,4,\dots ,w/num\_inputs\}$}
       \If{$(w/(\f{}*num\_inputs)) \bmod m_0\equiv 0$}
         \State $kernels\gets\emptyset$
         \ForAll{$i\in\{ 0,\dots ,num\_inputs-1\}$}
           \State insert $input\_groups[i][0:\f{}-1]$ in $kernels$
         \EndFor
         \State $s\gets\emptyset$
         \ForAll{$s\f{}\in kernels$}
           \State $s[s\f{}] = \text{VectorizedSumfactKernel}(kernels)$
         \EndFor
         \ForAll{$other\in$\Call{VectorizationStrategies}{$sum\f{}acts - kernels$, $width$, $q$}}
           \State $\text{extend}(strategies, \text{combine}(other, s))$
         \EndFor
       \EndIf
     \EndFor
   \EndFor
   \State\Return $strategies$
 \EndFunction
 \end{algorithmic}
\end{algorithm}

Algorithms~\ref{alg:vecopt} and \ref{alg:vecgen} establish a framework to explore different cost functions.
In our work, we have used two cost functions: An autotuning measurement function and a heuristic cost model function.
The autotuning function generates a benchmark program for a given sum factorization kernel, compiles it, and runs it on the target hardware.
Measured runtime is used as the cost function. While this approach delivers very good results and we generally use it for our performance measurements, there is also some downsides to it.
Depending on how large the set of possibilities is, code generation may take a substantial amount of time.
While this is not important for an HPC application to be run at large scale, it is unfeasible during a development cycle.
There is an additional problem with measuring sum factorization kernels as opposed to measuring whole integration kernels:
Cost increases in the quadrature loop are not covered by the benchmark programs.
To account for this, we leverage the fact that we can count floating point operations on the symbolic representation and introduce a penalty factor for the autotuning scores.

The heuristic cost model function we have worked with is given in equation~\eqref{equ:costfunction}.
It reproduces our practical experiences quite well, but does not take into account specific hardware features.
This can be used as a drop-in replacement, if the autotuning approach is infeasible.
The function depends on the following quantities:
\begin{itemize}
 \item the number of non-parallel floating point operations $\text{flops}(s\f{})$ carried out.
 \item a heuristic penalty function for the instruction level parallelism potential of a sum factorization kernel depending on the size of the splitting $s$ as used in section~\ref{sec:vertical}.
 where we observe that loop fusion based vectorization should always be preferred over vertical vectorization when applicable.
 \begin{equation}
  \text{ilp}(s\f{}) = 1 + c_0\log_2(s)
 \end{equation}

 \item a heuristic penalty function for the necessary load instructions. This depends on the number of input coefficients $p$ used for loop fusion in the kernel.
 \begin{equation}
  \text{loads}(s\f{}) = 1 + c_1\log_2(p)
 \end{equation}

\end{itemize}

\begin{figure}
 \includegraphics[width=0.3\textwidth]{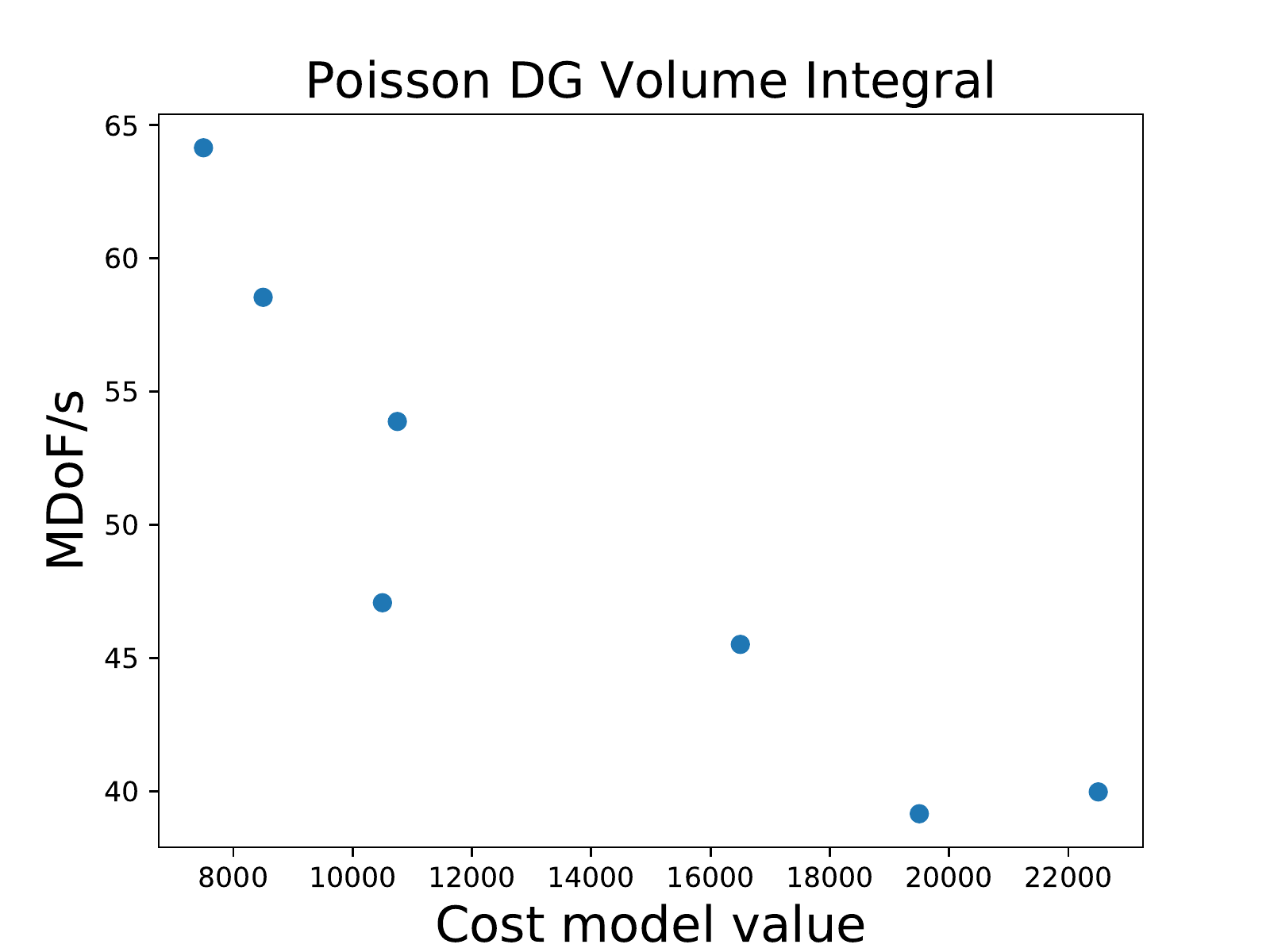}
 \includegraphics[width=0.3\textwidth]{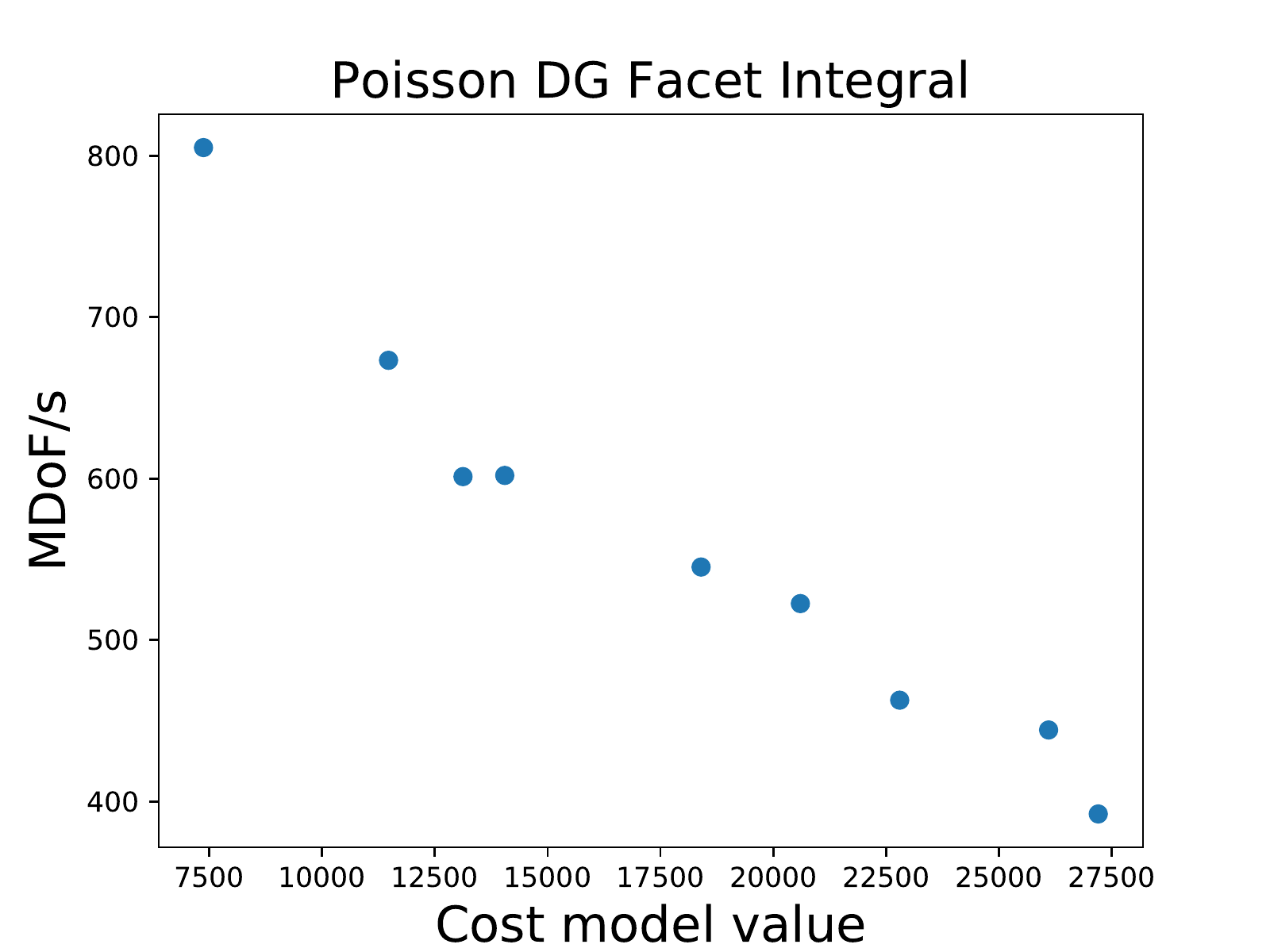}
 \includegraphics[width=0.3\textwidth]{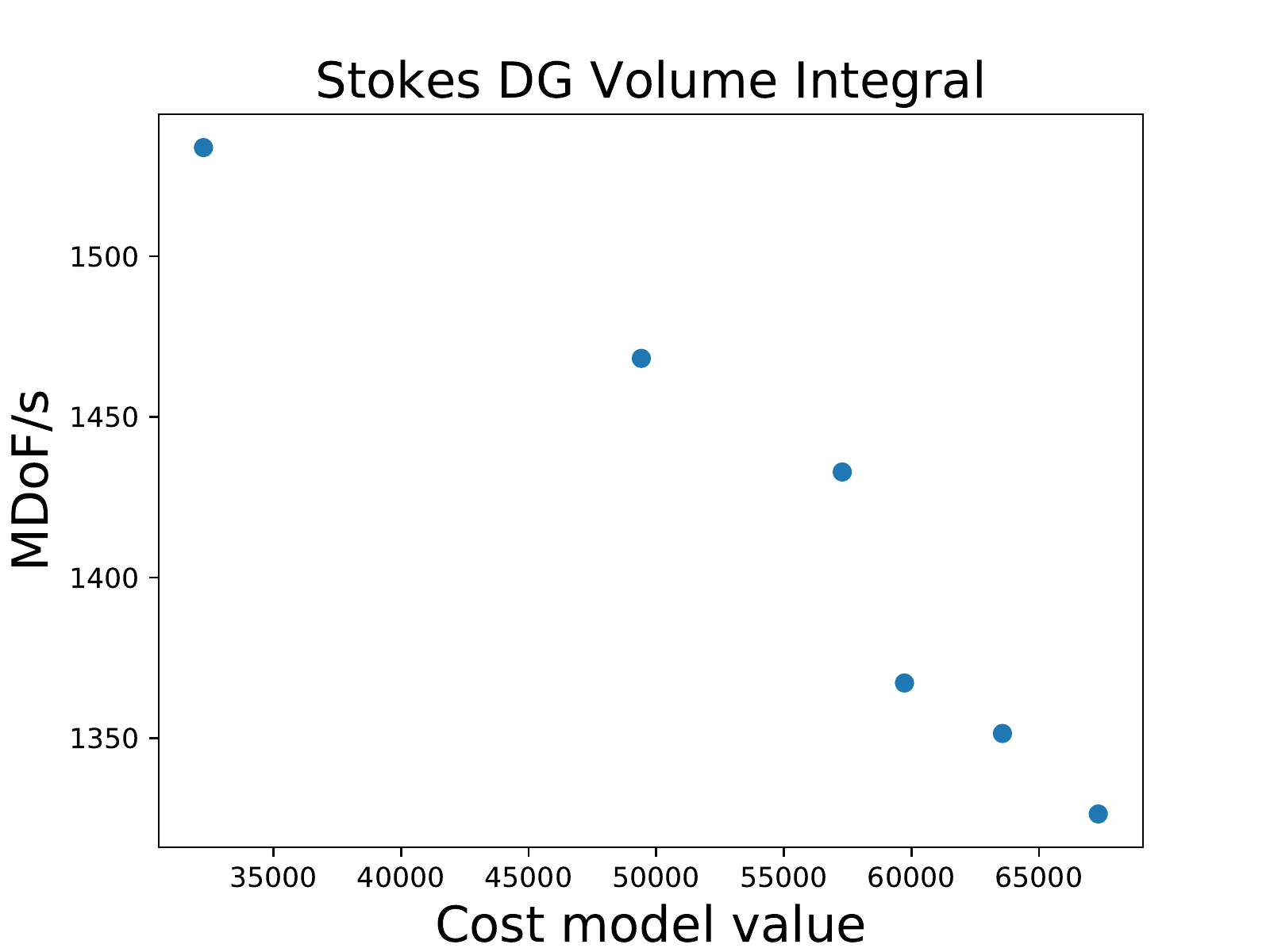}
 \caption{Validation of the cost model function from equation~\eqref{equ:costfunction}: Plotting the cost function value against measured performance.
 Each data point represents a possible vectorization strategy.
 The strategies are selected by running the optimization algorithm with a modified costfunction that minimizes the distance between the actual cost function value and a target value.}
 \label{fig:validation}
\end{figure}

The resulting cost function is the product of these given terms:
\begin{equation}
 \label{equ:costfunction}
 \text{cost}(s\f{}) = \text{flops}(s\f{}) \cdot \text{ilp}(s\f{}) \cdot \text{loads}(s\f{})
\end{equation}
In practice, we have chosen $c_0=0.5$ and $c_1=0.25$ and achieved good results.
Figure~\ref{fig:validation} validates this choice of cost model:
For a variety of integrals, a number of possible vectorization strategies is realized and performance is measured with the methodology described in section~\ref{sec:performance}.
We observe that the cost model minimum captures the performance maximum correctly.

\section{Performance Results}
\label{sec:performance}

Throughout this chapter, we will study the performance of the proposed algorithms on two different architectures.
As we strongly believe, that in order for performance numbers to be meaningful, all details about the used hardware and the measuring process must be given, we will thoroughly describe the benchmark setup in section~\ref{sec:benchsetup}.
We will then carry out performance studies for two model problems: Section~\ref{sec:diffreac} studies a diffusion-reaction type model, while section~\ref{sec:stokes} studies the Stokes equations.

\subsection{Benchmark Methodology}
\label{sec:benchsetup}

In this article we will use two measures to evaluate the performance of our implementation:
\begin{itemize}
 \item Floating Point Operations per second expressed in GFlops/s and often given as a percentage of the machine's maximum floating point performance.
 \item Degrees of freedom per second processed during a full operator application. Note that we prefer this measure over its inverse (time per degree of freedom).
\end{itemize}

Good results on the latter measure are always more important from the application point of view, as it gives an accurate measure of how fast a real problem can be solved.
However, the former is still an interesting measure that allows reasoning about how good a code is suited for a given hardware. \\

In order to accurately measure the floating point operations per second, the number of performed floating point operations needs to be measured exactly.
We exploit the fact that we are generating code for dune-pdelab, which uses C++ templates to the extent that we can replace the underlying floating point type throughout all our simulation workflow.
Instead of using double, we use a custom type templated to double, which has overloads for all arithmetic operations that increase a global counter and forward the operation to the underlying template type.
This counting of course introduces a non-negligible performance overhead. We therefore compile different executables from the same source for operation counting and time measurement. \\

Apart from counting operations, accurate time measurements are needed.
We instrument our code with C macros to start and stop timers using the TSC registers.
The performance overhead of starting and stopping a timer is measured at runtime and the measurements are calibrated accordingly.
To gain further insight into the performance bottlenecks of our implementation, we measure time and floprates at different levels of granularity:
Full operator application, cell-local integration kernel and individual steps of algorithm~\ref{alg:assembly}.
For all these granularity levels, separate executables are compiled to assure that no measurement is tampered by additional measurements taken within the measuring interval. \\

We study the node-level performance of our generated code on two Intel micro architectures.
\begin{itemize}
 \item Intel Haswell
 \begin{itemize}
 \item Intel Xeon processor E5-2698v3
 \item 2x16 cores
 \item base frequency: 2.3 GHz, 1.9 Ghz on AVX2-heavy loads \cite{dolbeau2017theoretical}
 \item Theoretical peak performance: 972.8 GFlops/s0
 \end{itemize}
 \item Intel Skylake
 \begin{itemize}
 \item Intel Xeon Gold 6148
 \item 2x20 cores
 \item base frequency: 2.4 GHz, 2.0 GHz on AVX2-heavy loads, 1.6 GHz on AVX512-heavy loads \cite{dolbeau2017theoretical}
 \item Theoretical peak performance: 2.04 TFlops/s
 \end{itemize}
\end{itemize}

We turn off turbo mode entirely on the machines in order to be able to get a good estimate on the maximum floating point performance of the machine.
In order to study the node-level performance of our implementation, we need to saturate the entire node with our computation.
Otherwise, some processors may have priviliged access to ressources such as memory controllers and tamper results.
We do so by doing MPI parallel computations with one rank per processor and an identical workload size on each of these ranks.
Also, we choose the problem size to be such, that one vector of degrees of freedom exceeds the level 3 cache of the machine.
While this may not be a realistic setting when doing strong scaling of simulation codes, it gives a good worst case estimate of our code's node level performance.
The time for communication of overlap data via MPI is not included in our measurement, as the task of hiding that communication behind computation is not the subject of this work.

\subsection{Diffusion-Reaction Equation}
\label{sec:diffreac}
In this section, we will consider the diffusion-reaction equation~\eqref{equ:diffreac} on $\Omega =(0,1)^3$, where we assume the triangulation $\mathcal{T}_h$ to be axiparallel.

\begin{align}
\label{equ:diffreac}
 -\nabla\cdot\left( \mathbf{k}(\mathbf{x})\nabla u(\mathbf{x})\right) + c(\mathbf{x})u(\mathbf{x}) &= f(\mathbf{x}) && \forall\mathbf{x}\in\Omega \\
 u(\mathbf{x}) &= g(\mathbf{x}) && \forall\mathbf{x}\in\partial\Omega\nonumber
\end{align}

For the discretization we choose the symmetric interior penalty Discontinuous Galerkin method, that we recall in equation~\eqref{equ:diffreac_weak} to disclose all the details of our implemenation to enable the reader to compare our performance numbers to other codes.

First we need to introduce some notation.
For an interior facet $F$ with neighboring elements $T_1$ and $T_2$ we define the jump
\begin{equation*}
  \jump{v}(x) := v|_{T_1}(x) - v|_{T_2}(x)
\end{equation*}
the average
\begin{equation*}
  \avg{v}(x) := \frac{1}{2}\left(v|_{T_1}(x) + v|_{T_2}(x)\right).
\end{equation*}
and the penalty factor
\begin{equation*}
  \gamma_F := \frac{\alpha k(k+d-1) |F|}{\min (|T^+(F)|, |T^-(F)|)}
\end{equation*}
where $k$ is the polynomial degree (same for all directions), $d$ is the world dimension and $\alpha$ a free parameter we set to $3$. On a boundary facet $F$ belonging to an element $T$ we define
\begin{equation*}
  \jump{v}(x) := \avg{v}(x) := v|_{T}(x).
\end{equation*}
and
\begin{equation*}
  \gamma_F = \frac{\alpha k(k+d-1) |F|}{|T|}
\end{equation*}
respectively. Then the residual function corresponding to problem~\eqref{equ:diffreac} is
\begin{align}
  r_h(u_h,v_h)
  &= \sum_{T\in\mathcal{T}_h} \left[
    (\mathbf{k} \nabla u_h, \nabla v_h)_{0,T}
    + (cu_h,v_h)_{0,T}
    - (f,v_h)_{0,T}
    \right] \nonumber \\
  &\quad + \sum_{F\in\mathcal{F}_h^i \cup \mathcal{F}_h^b} \left[
    -(\mathbf{n}_F \cdot \avg{\mathbf{k} \nabla u_h}, \jump{v_h})_{0,F}
    + \theta (\jump{u_h},\mathbf{n}_F \cdot \avg{\mathbf{k} \nabla v_h})_{0,F}
    + \gamma_F (\jump{u_h}, \jump{v_h})_{0,F}
    \right] \nonumber \\
  &\quad + \sum_{F\in\mathcal{F}_h^b} \left[
    - \theta (g, \mathbf{n}_F \cdot (\mathbf{k} \nabla v_h))_{0,F}
    - \gamma_F (g, v_h)_{0,F} \right],  \label{equ:diffreac_weak}
\end{align}
where we choose $\theta=-1$.

We choose the physical parameter functions in equation~\eqref{equ:diffreac} in such a way that the function $g$ is a manufactured solution:
\begin{align}
 \mathbf{k}(\mathbf{x}) & = \mathbf{x}\mathbf{x}^T + \mathbf{I} \\
 c(\mathbf{x}) & = 10 \\
 f(\mathbf{x}) & = -6 \\
 g(\mathbf{x}) & = \|\mathbf{x}\|_2^2
\end{align}

\begin{figure}
 \includegraphics[width=0.49\textwidth]{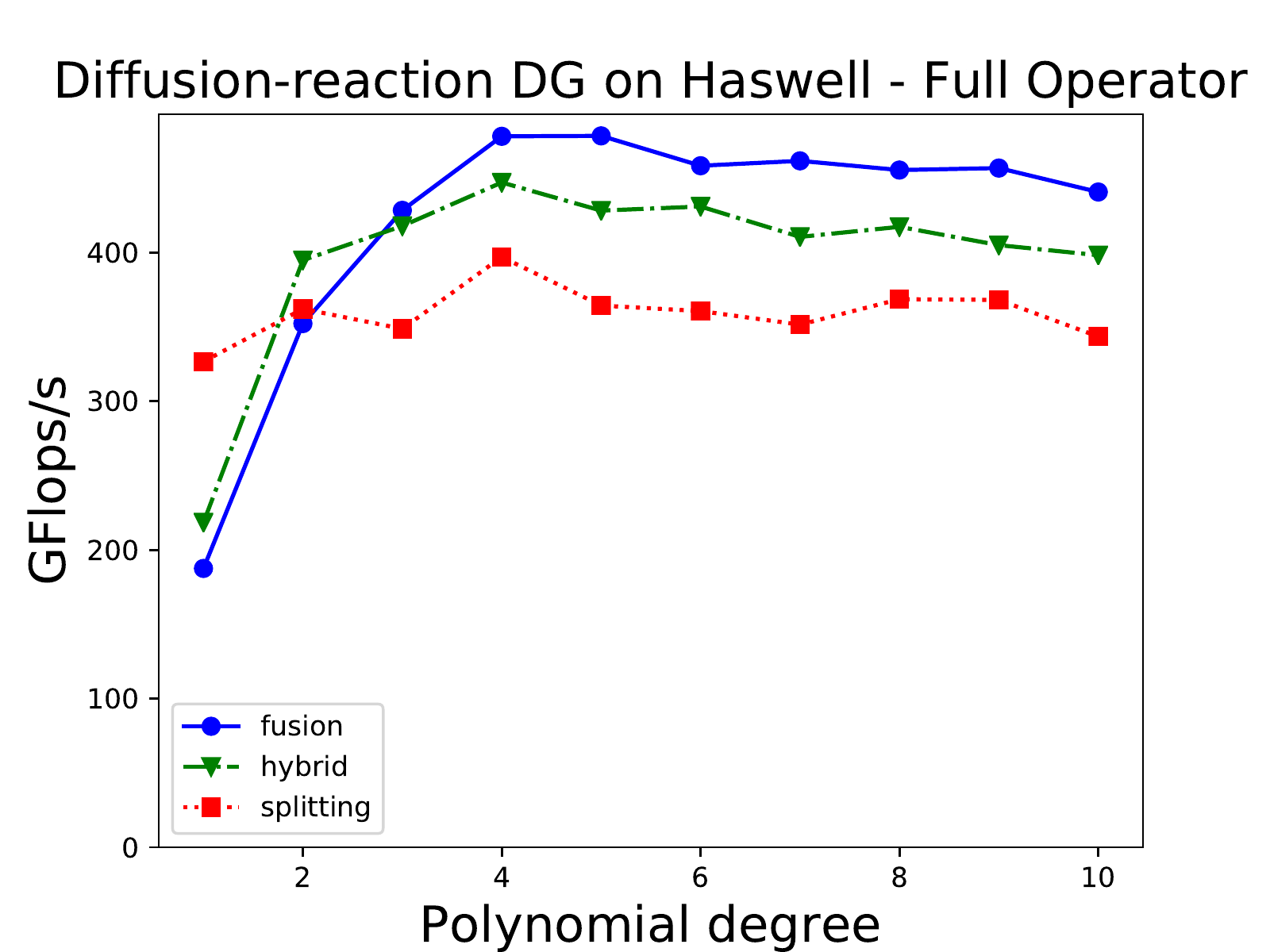}
 \includegraphics[width=0.49\textwidth]{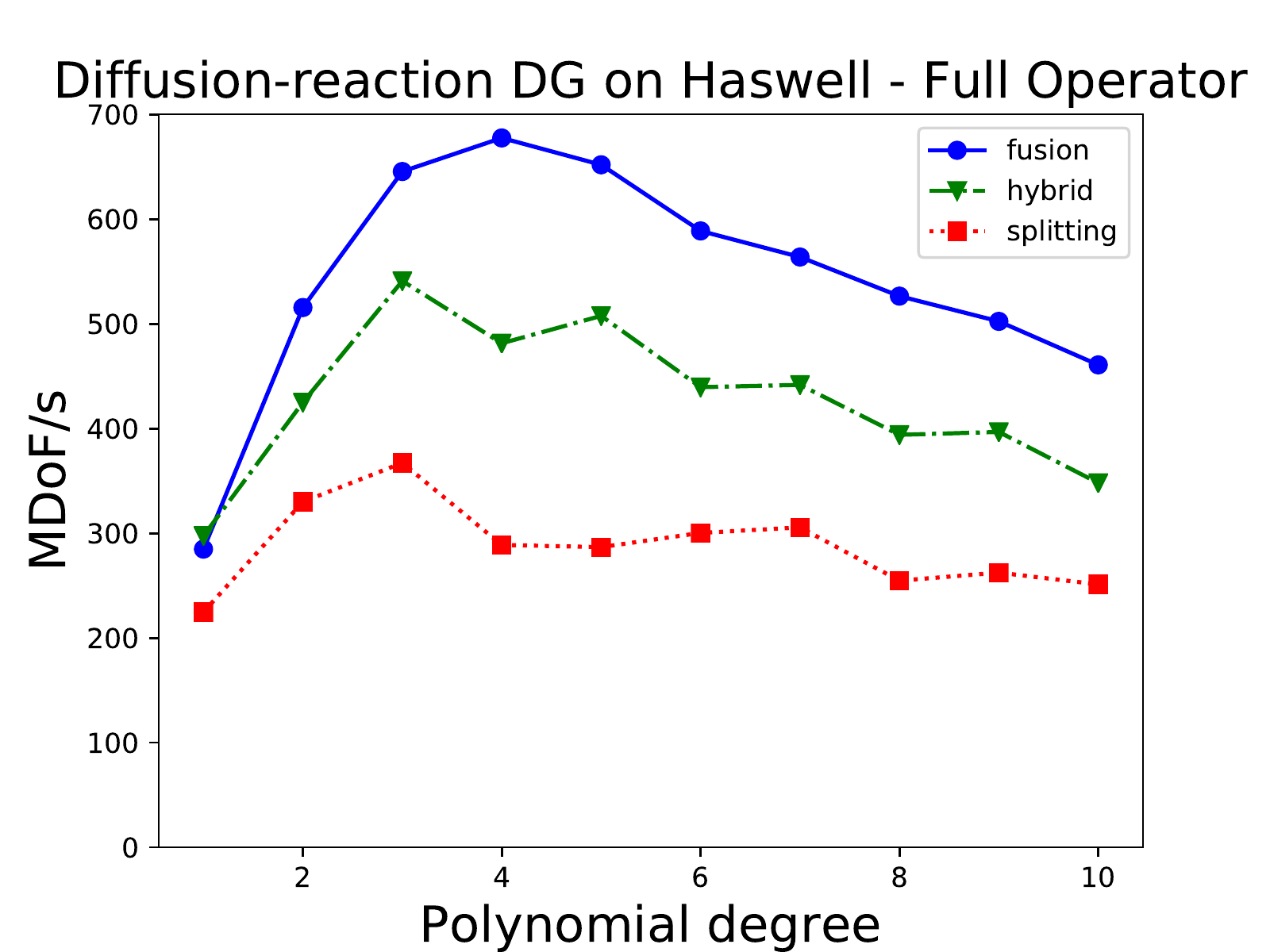}
 \caption{Performance numbers for the operator application of the Poisson DG problem described in equation~\eqref{equ:diffreac_weak} on the Intel Haswell system described in section~\ref{sec:benchsetup}.
 Left plot shows the GFlops/s achieved, right plot the degrees of freedom per second, both for a full operator application. The employed vectorization strategies are taken
 from section~\ref{sec:vectorization} with the hybrid strategy combining two sum factorization kernels into a 256 bit SIMD register.}
 \label{fig:haswell_poisson}
\end{figure}

Figure~\ref{fig:haswell_poisson} shows the performance results for an operator application of equation~\eqref{equ:diffreac_weak} on the Haswell system described in section~\ref{sec:benchsetup} with a variety of vectorization strategies from section~\ref{sec:vectorization} being applied.
We observe that - given a sufficiently high polynomial degree to hit the flopbound regime - our approach is able to get sustained 50\% of machine peak on the given Intel Haswell node.
Given the fact that we are measuring a full operator application and not only the sum factorization kernels, this is a highly competitive value.
Comparing the described vectorization strategies, we note that loop fusion based vectorization clearly outperforms the other strategies once we are in the flopbound regime.
This is most likely due to the fact, that the size of the loop bounds of the vectorized kernels is larger by a factor of two and four respectively compared to hybrid and splitting based strategies.
This increased size allows superior instruction level parallelism in terms of unrolling and data parallelism between the two floating point units of the core.
We conclude, that loop fusion based vectorization should be the strategy of choice whenever the given PDE problem allows it by requiring a suitable number of quantities to be evaluated.
This conclusion has already been incorporated into the cost model approach in section~\ref{sec:sumfactcodegen}. \\

\begin{figure}
 \includegraphics[width=0.49\textwidth]{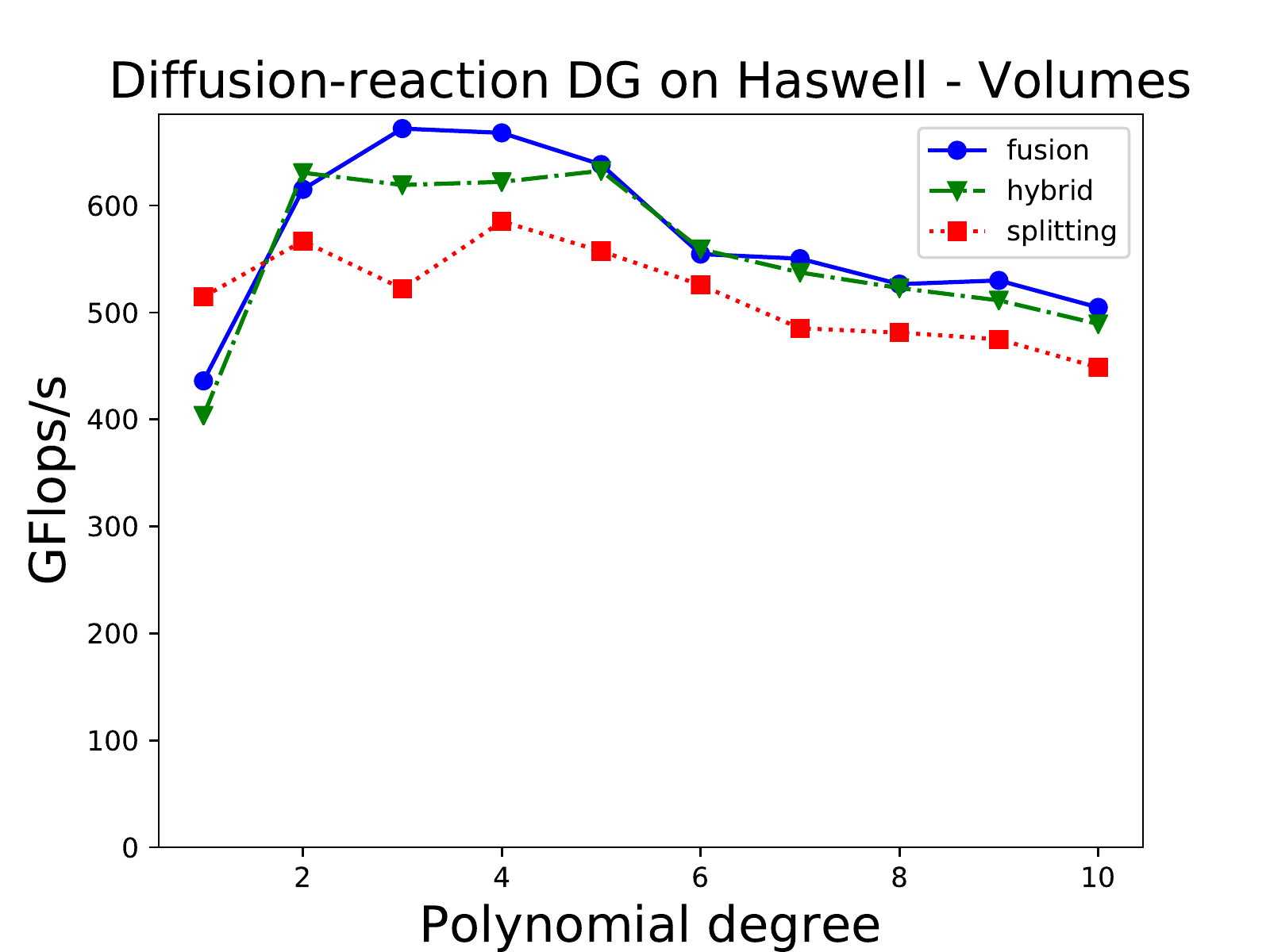}
 \includegraphics[width=0.49\textwidth]{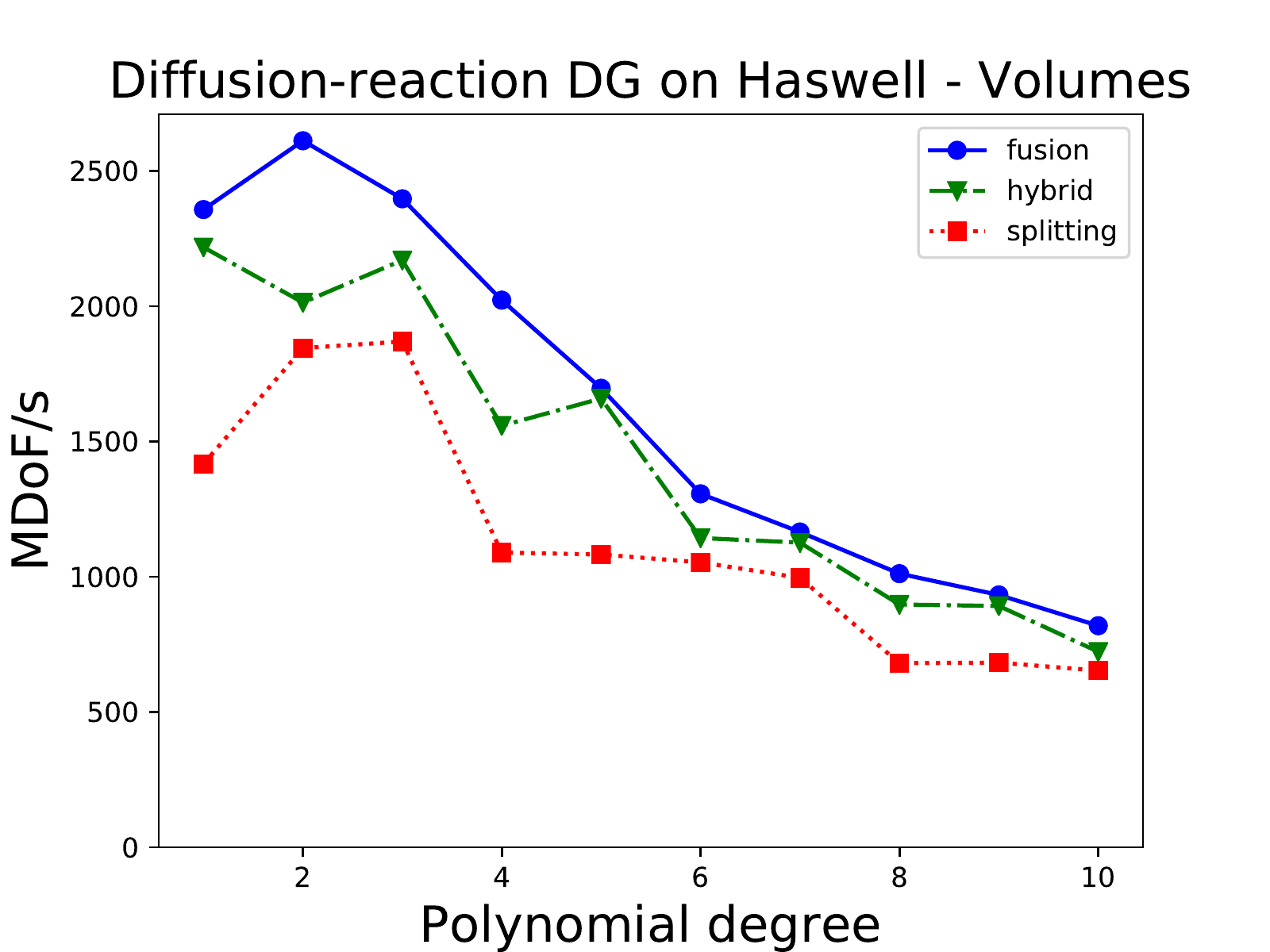}
 \caption{Performance numbers for the volume integrals of an operator application of the Poisson DG problem described in equation~\eqref{equ:diffreac_weak} on the Intel Haswell system described in section~\ref{sec:benchsetup}.}
 \label{fig:haswell_poisson_volume}
\end{figure}

\begin{figure}
 \includegraphics[width=0.49\textwidth]{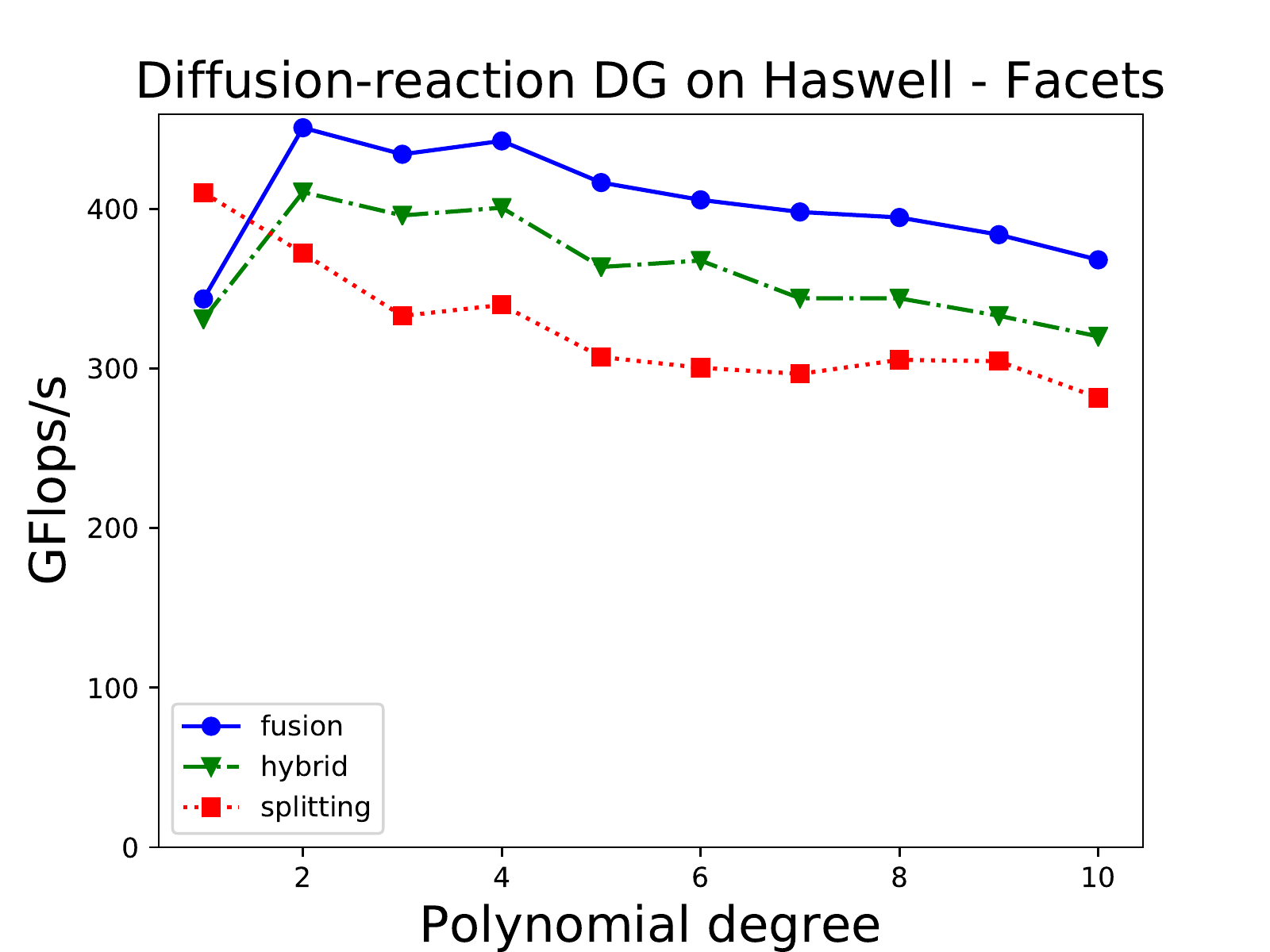}
 \includegraphics[width=0.49\textwidth]{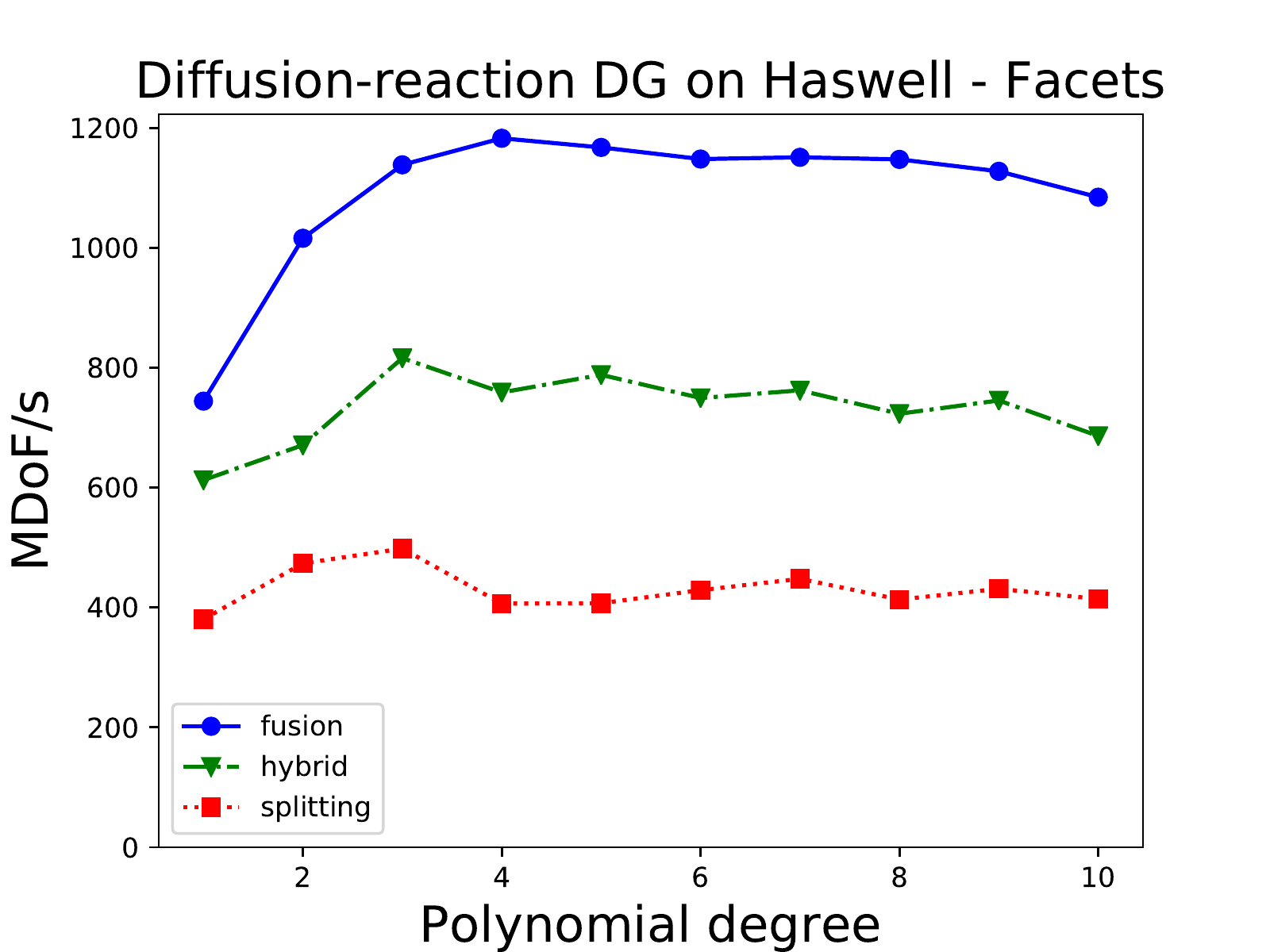}
 \caption{Performance numbers for the facet integrals of an operator application of the Poisson DG problem described in equation~\eqref{equ:diffreac_weak} on the Intel Haswell system described in section~\ref{sec:benchsetup}.}
 \label{fig:haswell_poisson_skeleton}
\end{figure}

Looking into the performance of our implementation at a more fine grained level, we can study differences between volume and facet integral contributions to the operator application.
We do not study boundary integrals as they are very similar to the integrals on interior facets and their total contribution is negligible.
Figure~\ref{fig:haswell_poisson_volume} shows the performance numbers for volume integrals, figure~\ref{fig:haswell_poisson_skeleton} for facet integrals.
We note that due to their favorable ratio of floating point operations to memory accesses, volume integrals achieve a signifcantly higher percentage of machine peak.
Also, vertical vectorization suffers from a higher performance loss for skeleton integrals.
This can again be explained with the size of the loop nests and the resulting instruction level parallelism capabilities. \\

\begin{figure}
 \includegraphics[width=0.49\textwidth]{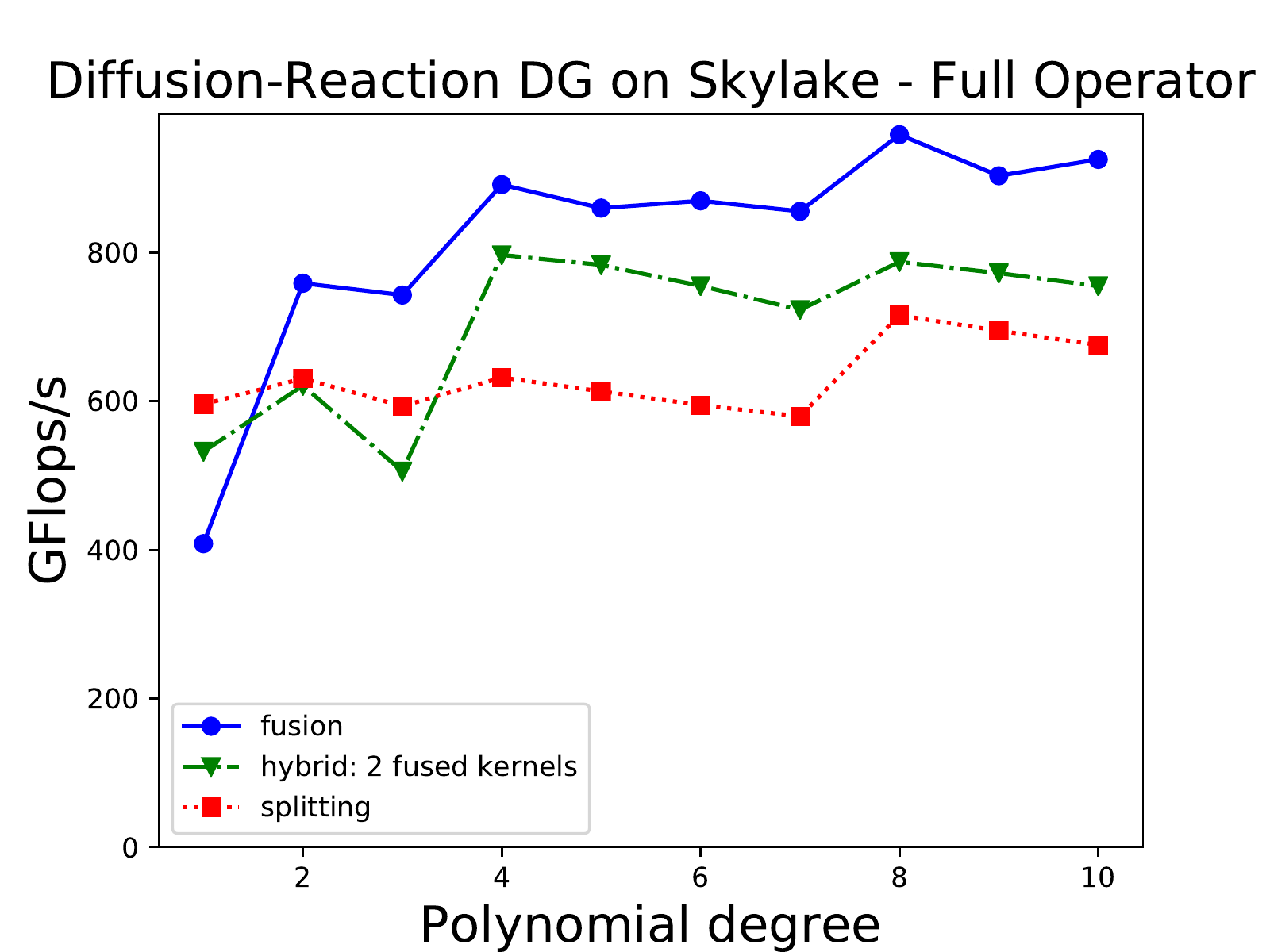}
 \includegraphics[width=0.49\textwidth]{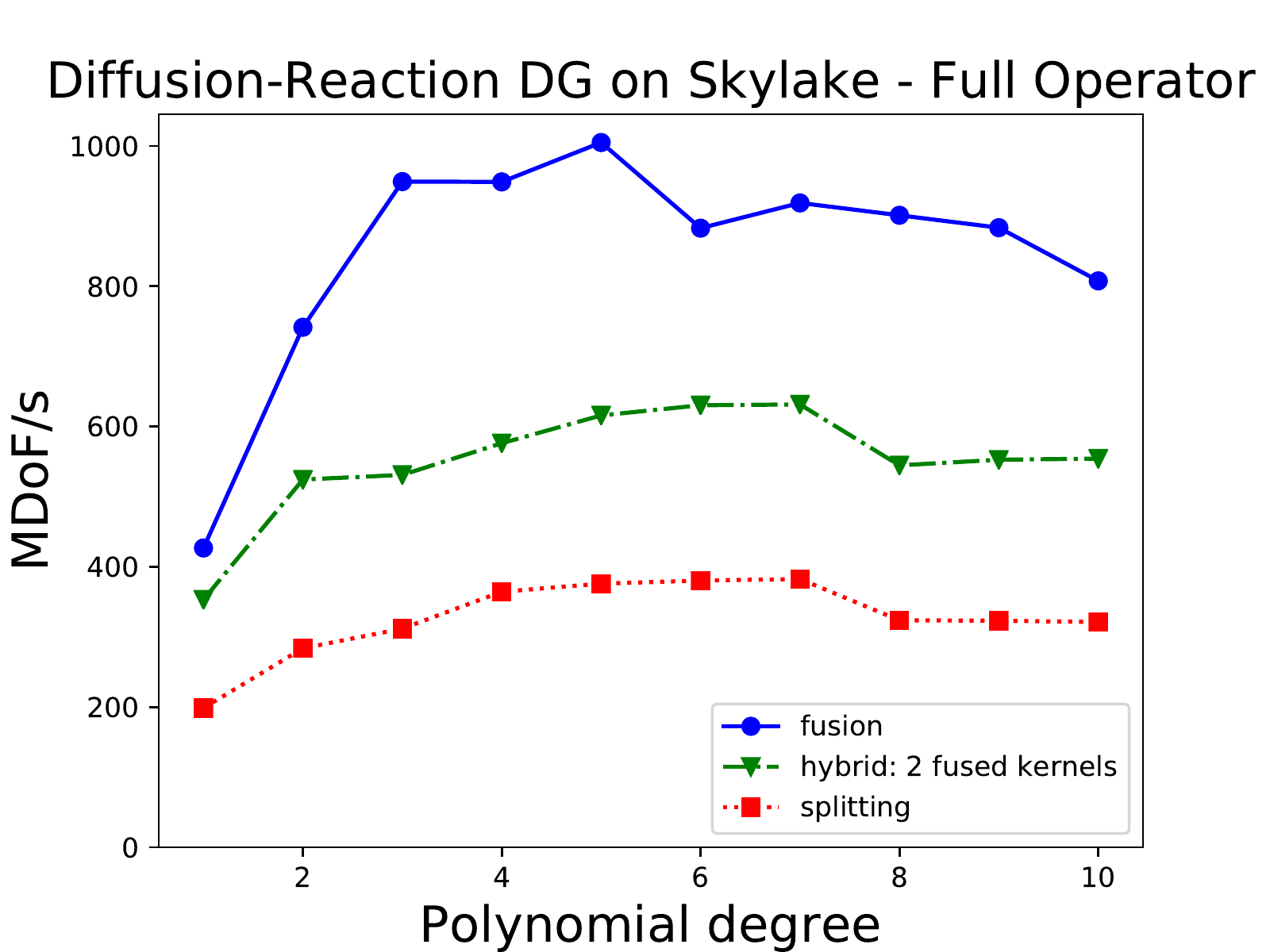}
 \caption{Performance numbers for the operator application of the Poisson DG problem described in equation~\eqref{equ:diffreac_weak} on the Intel Skylake system described in section~\ref{sec:benchsetup}.
          Results for a variety of strategies using loop fusion and loop splitting.
 }
 \label{fig:skylake}
\end{figure}

\begin{figure}
 \includegraphics[width=0.49\textwidth]{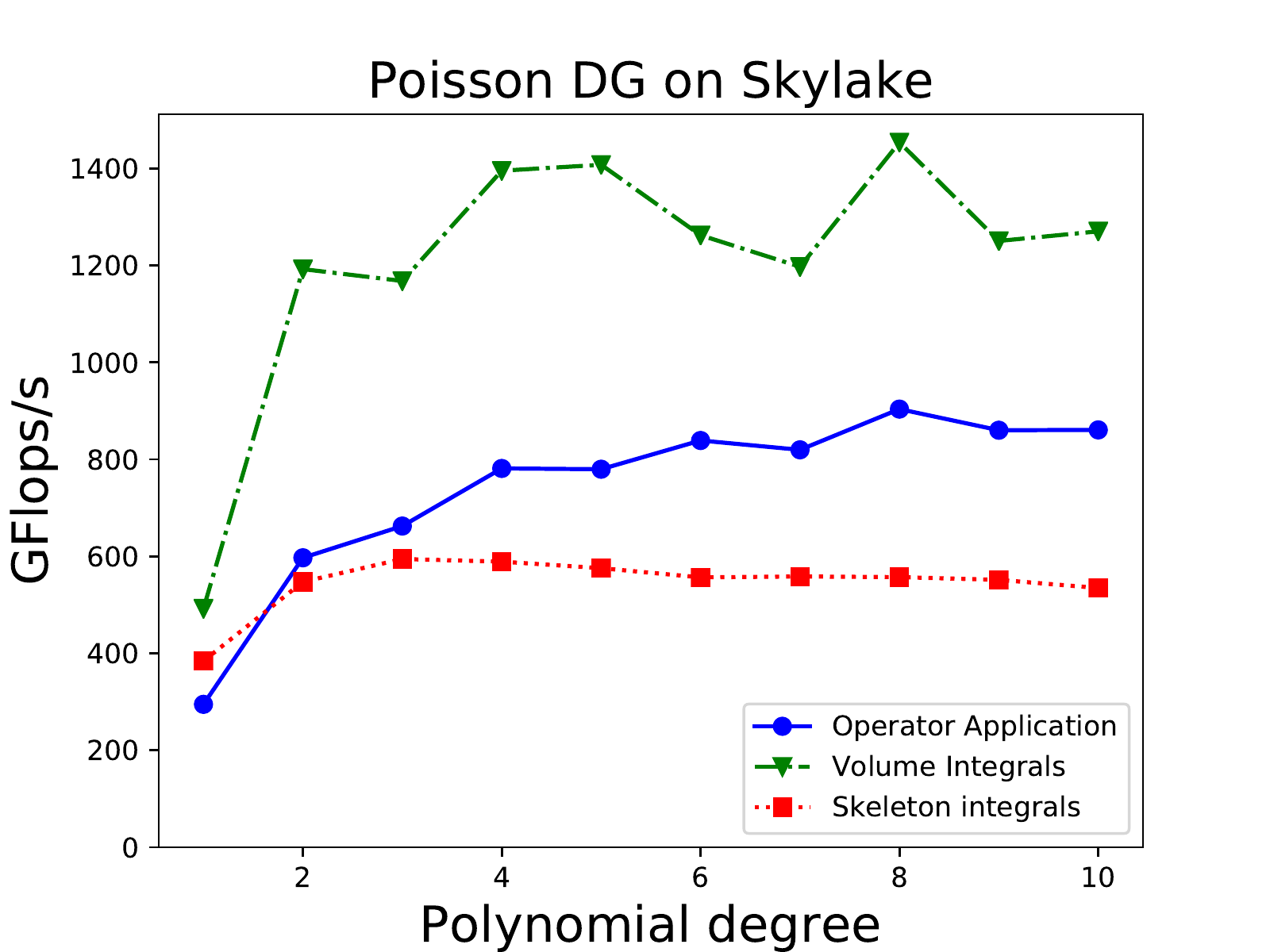}
 \includegraphics[width=0.49\textwidth]{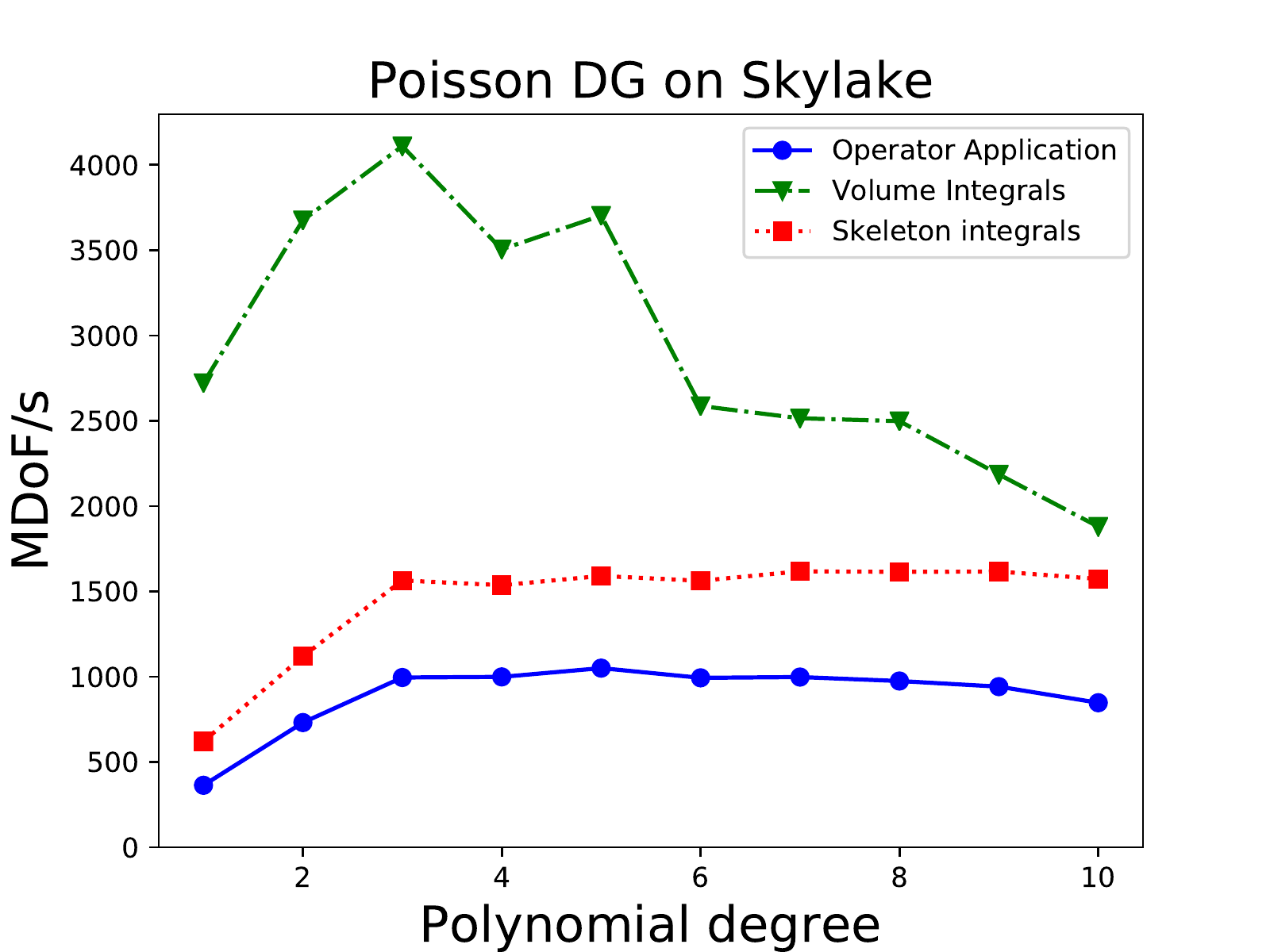}
 \caption{Performance numbers for the operator application of the Poisson DG problem described in equation~\eqref{equ:diffreac_weak} on the Intel Skylake system described in section~\ref{sec:benchsetup}.
          Results with the costmodel from equation~\eqref{equ:costfunction}, given for volume and facet integrals individually.}
 \label{fig:skylake2}
\end{figure}

Figure~\ref{fig:skylake} shows the performance numbers for the Poisson problem on the given Intel Skylake system.
Again, we first have a look at how fusion based strategies compare to splitting based ones.
We again come to the conclusion that fusion based approaches are favorable wherever applicable.
Having a look at absolute numbers, we note that for medium and high polynomial degrees we get >40\% of the machine peak for an entire operator application.
As this is a slightly lower percentage than we were achieving on Haswell, we will look a bit into detail here:
Figure~\ref{fig:skylake2} shows the vectorization strategy picked by the costmodel split by volume and facet integrals.
We see that the volume integrals are capable of exploiting the hardware really well - with a percentage of up to 70\%.
However, the facet integrals are only capable of exploting roughly 25\% of machine peak.
There are several reasons for this inferior performance, where some of these are caused by the nature of facet integrals and some are caused by our implementation:
\begin{itemize}
 \item The loop bounds of sum factorization kernels on facets are of lower dimensionality - disallowing a lot of data parallelism between the floating point units of the core.
 \item The vectorization strategy fuses the evaluations of $u^+$, $\nabla u^+$, $u^-$, $\nabla u^-$, which requires an 8x8 matrix transposition analoguous to figure~\ref{fig:transpose}.
       Such a transposition is much more costly than the 4x4 one and may be harder to hide behind the quadrature loop floating point operations.
\end{itemize}
Despite this improvable performance for facet integrals, we are still able to get a better throughput in terms of degrees of freedoms per second from the Skylake node, than we get from the Haswell node.

\subsection{Stokes Equation}
\label{sec:stokes}

As a second example we consider the steady Stokes equation of incompressible viscous flows on the domain $\Omega=(0,1)^3$ and apply the vectorization strategies of chapter~\ref{sec:vectorization}.
The strong formulation of this problem is given by
\begin{align}
  \label{equ:stokes_strong}
  - \mu \Delta \mathbf{u} + \nabla p &= \mathbf{f} \qquad \text{in } \Omega \\
  \nabla \cdot \mathbf{u} &= 0 \qquad \text{in } \Omega \\
  \mathbf{u} &= \mathbf{g} \qquad \text{on } \Gamma_D \\
  \mu \nabla \mathbf u \mathbf{n} - p \mathbf{n} &= 0 \qquad \text{on } \Gamma_N
\end{align}
for Dirichlet boundary $\Gamma_D\subset\partial\Omega$, $\Gamma_D\neq\{\emptyset, \partial\Omega\}$ and Neumann boundary $\Gamma_N=\partial\Omega \setminus \Gamma_D$.
The residual function of the Discontinuous Galerkin discretization of this problem is
\begin{align}
  r(\mathbf{u}_h, p_h, \mathbf{v}_h, q_h)
  &= \sum_{T\in\mathcal{T}_h} \left[
    (\nabla \mathbf{u}_h, \nabla\mathbf{v}_h)_{0,T}
    - (p_h, \nabla \cdot \mathbf{v}_h)_{0,T}
    - (\nabla \cdot \mathbf{u}_h, q)_{0,T}
    \right] \nonumber \\
  &\quad + \sum_{F\in\mathcal{F}_h^i \cup \mathcal{F}_h^D} \left[
    - (\avg{\nabla \mathbf{u}_h} \mathbf{n}_F, \jump{\mathbf{v}_h})_{0,F}
    + \theta (\jump{\mathbf{u_h}}, \avg{\mathbf{v}_h} \mathbf{n}_F)_{0,F}
    + \gamma_F (\jump{\mathbf{u}_h},\jump{\mathbf{v}_h})_{0,F}
    \right] \nonumber \\
  &\quad + \sum_{F\in\mathcal{F}_h^i \cup \mathcal{F}_h^D} \left[
    (\avg{p_h}, \jump{\mathbf{v}_h} \mathbf{n}_F)_{0,F}
    + (\jump{\mathbf{u}_h} \mathbf{n}_F, \avg{q_h})_{0,F}
    \right] \nonumber \\
  &\quad + \sum_{F\in\mathcal{F}_h^D} \left[
    - \theta (\mathbf{g}, \nabla \mathbf{v}_h \mathbf{n}_F)_{0,F}
    - \gamma_F (\mathbf{g}, \mathbf{v}_h)_{0,F}
    - (\mathbf{g} \cdot \mathbf{n}_F, q_h)_{0,F}
    \right], \label{equ:stokes_weak}
\end{align}
where $\mathcal{F}_h^D=\mathcal{F}_h^b\cap\Gamma_D$ is the set of boundary faces with Dirichlet boundary condition.
We assume that each boundary face is either completely part of the Dirichlet or Neuman boundary.
We choose the parameters $\theta=-1$ and $\gamma_F$ as in the previous section.
This is a system of equations with solution velocity $\mathbf{u}_h\in (V_h^k)^3$ and pressure $p_h\in V_h^{k-1}$ and test functions $\mathbf{v}_h\in (V_h^k)^3$ and $q_h\in V_h^{k-1}$ where the polynomial degree $k$ is again the same for each direction.

For the calculations below we use $\mu=1$, $\Gamma_D=\{x\in\partial\Omega:x_1<1\}$ and $\mathbf{g}=(4x_2(1-x_2),0,0)^T$. This way $\mathbf{g}$ is again a manufactured solution for the velocity and the solution for the pressure is given by $p=8(1-x_1)$.



\begin{figure}
 \includegraphics[width=0.49\textwidth]{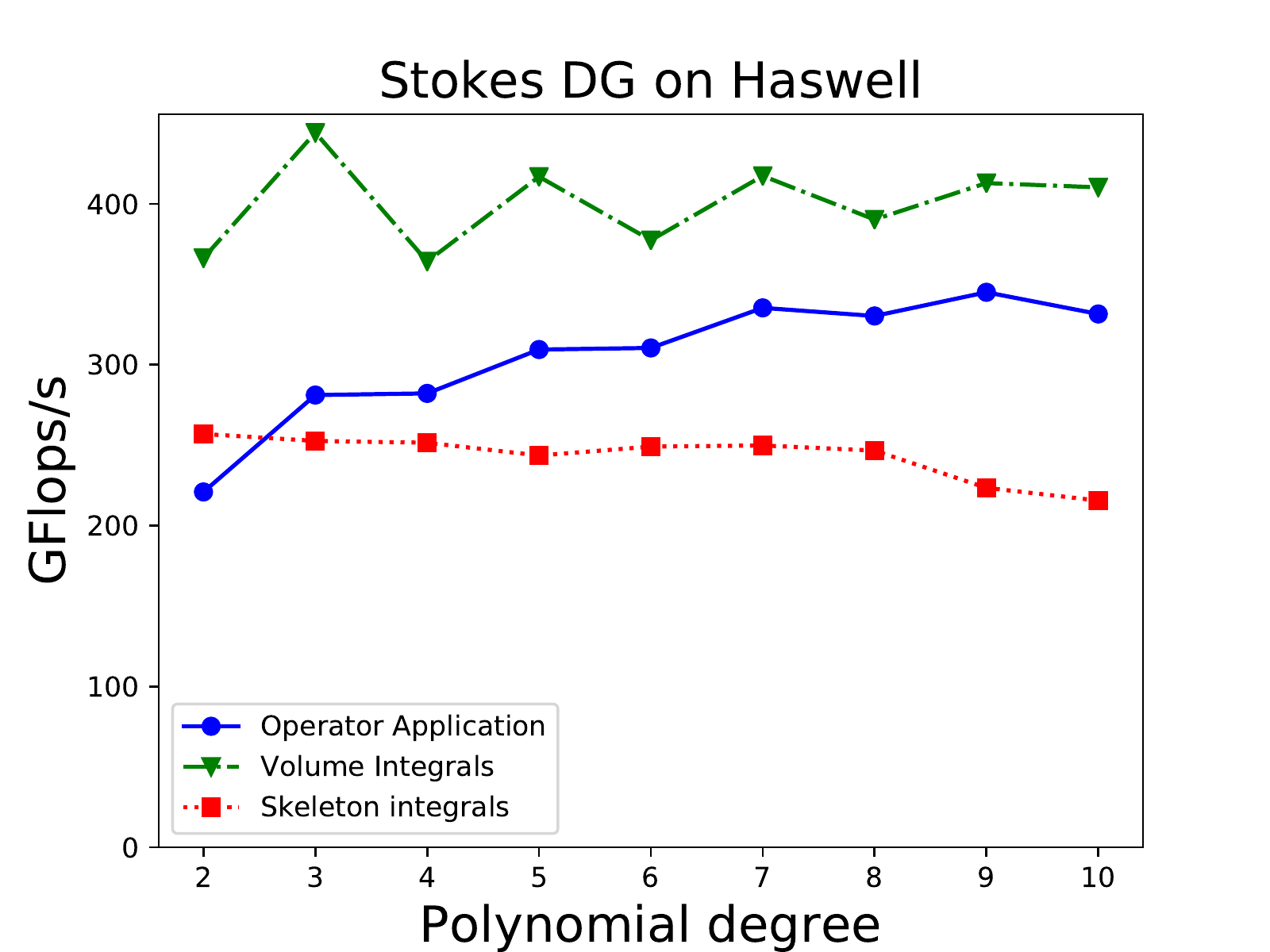}
 \includegraphics[width=0.49\textwidth]{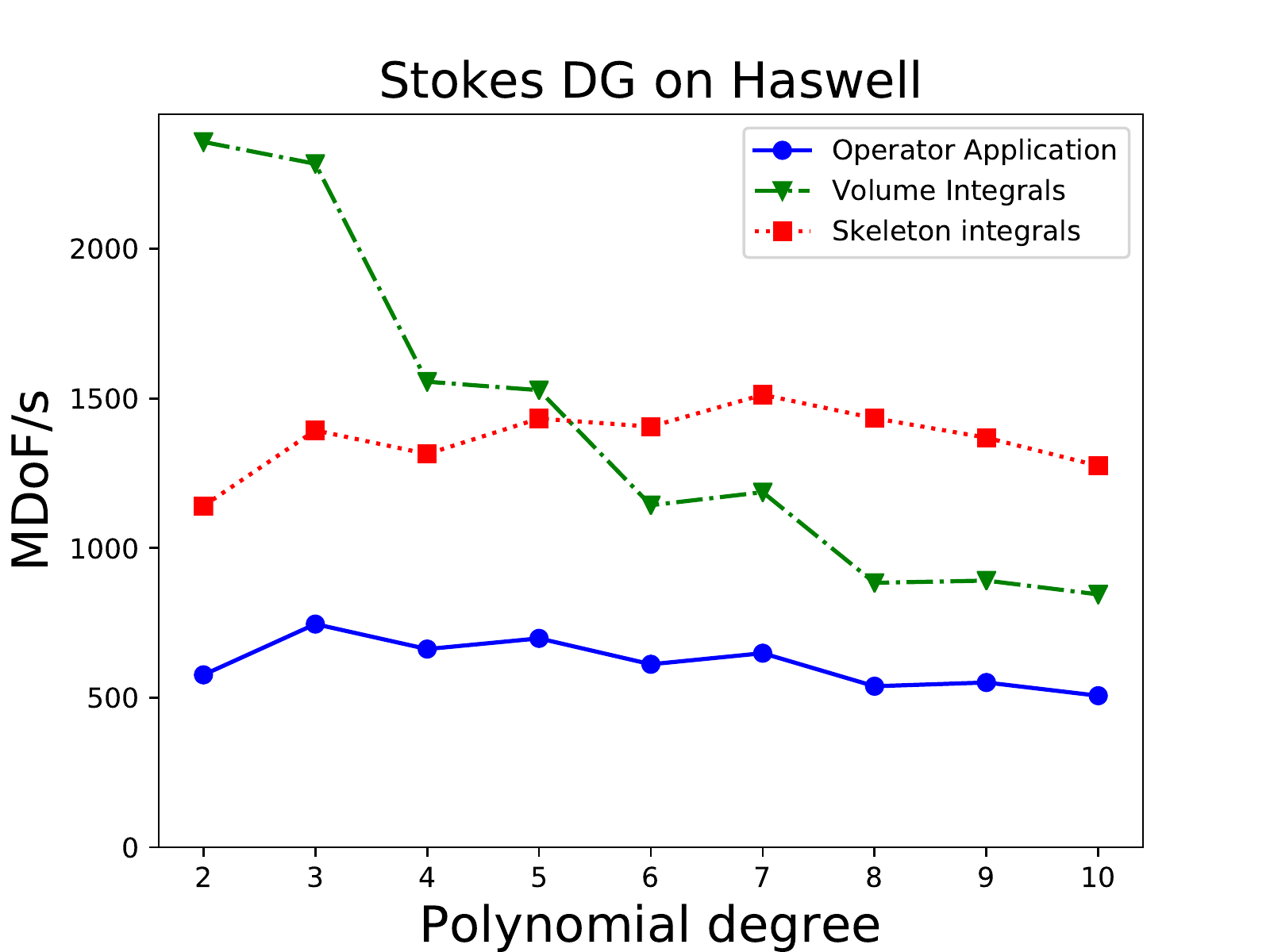}
 \caption{Performance numbers for the Stokes DG problem described in equation~\eqref{equ:stokes_weak} on the Intel Haswell system described in section~\ref{sec:benchsetup}.
 Vectorization is done through the cost model approach described in section~\ref{sec:sumfactcodegen}.
 The $\mathcal{Q}_2/\mathcal{Q}_1$ data point on the left hand side suffers from a measurement artifact: The more fine grained measurements do not sum correctly to the coarse grained measurements resulting the flop rate for the overall application not being between the volume and facet flop rate.
 For higher polynomial degrees, this effect vanishes.
 }
 \label{fig:haswell_stokes}
\end{figure}

We show the performance of our implementation of the 3D Stokes problem given in equation~\eqref{equ:stokes_weak} in figure~\ref{fig:haswell_stokes}.
Vectorization is done using the cost model approach described in section~\ref{sec:sumfactcodegen}.
The cost model always picks a fusion based vectorization strategy for the evaluation of $\nabla\mathbf{u}$ without raising the number of quadrature points over the user defined minimum input on volume integrals.
On facets - with the grid being axiparallel - $(\nabla \mathbf{u}\mathbf{n})_i$ simplifies to $\partial_du_i$ with $d$ being the direction of the facet normal.
This results in a fusion based vectorization on facets with different input tensors in the lower and upper SIMD register half: $\partial_du^+|u^+|\partial_du^-|u^-$.
For the evaluation of pressure, the implementation depends on the divisibility of $m_0$: It is splitting based for $k=3$ and $k=7$ for volume integrals and scalar in for other polynomial degrees.
On facet integrals, a hybrid approach fusing the evaluation of $p^+$ and $p^-$ is applied for all odd polynomial degrees.
These cost model decisions are clearly reflected in the performance numbers in equation~\ref{fig:haswell_stokes}.

Comparing the GFlop rates of the Stokes problem to those of the Poisson problem, we have to use the hybrid strategy from figure~\ref{fig:haswell_poisson}, as the baseline for comparison.
Comparison against figure~\ref{fig:haswell_poisson} shows a reduced, but not drastically reduced, floating point throughput of approximately 30\% of machine peak.
This can be explained with the increased working set size for a PDE system and the resulting cache effects. \\

\section{Conclusion and future work}

We have presented a class of algorithms based on loop fusion and loop splitting that allow explicit vectorization of the finite element assembly problem with high order Discontinuous Galerkin methods.
The presented results show that we can utilize a notable percentage of the machines floating point capabilities for the Intel Haswell and Skylake architectures.
The algorithms implementation is embedded into a code generation workflow, that allows very flexible implementations depending on the given PDE and the given hardware.
Decisions on vectorization are delegated to a costmodel approach.
The authors will continue working on various aspects of this workflow.
The cost model function from equation~\eqref{equ:costfunction} can be played around with and extended to achieve even better vectorization.
The intermediate representation loopy allows for a lot more experiments to optimize the assembly loop nests.
Implementing facet integrals with sum factorization on unstructured grids is currently implemented into our codebase.
Future applications will include a massively parallel solver for subsurface flow.

\section{Acknowledgments}

The authors would like to acknowledge the impressive work of Andreas Kl\"ockner on the loopy package and its depedencies.
This work was partially funded by the German Federal Ministry of Education and Research (BMBF) under grant 01IH16003C.
Responsibility for the content of this publications lies with the authors.

\bibliographystyle{plain}
\bibliography{./perftool-bib/perftool.bib,./perftool-bib/fenics.bib}

\end{document}